# Accurate calculations of Stationary Distributions and Mean First Passage Times in Markov Renewal Processes and Markov chains

Jeffrey J. Hunter

Department of Mathematical Sciences
School of Engineering, Computer and Mathematical Sciences
Auckland University of Technology, New Zealand
Email: jeffrey.hunter@aut.ac.nz

**Abstract**
This article describes an accurate procedure for computing the mean first passage times of a finite irreducible Markov chain and a Markov renewal process. The method is a refinement to the Kohlas, *Zeit fur Oper Res*, 30,197-207, (1986) procedure. The technique is numerically stable in that it doesn't involve subtractions. Algebraic expressions for the special cases of one, two, three and four states are derived. A consequence of the procedure is that the stationary distribution of the embedded Markov chain does not need to be derived in advance but can be found accurately from the derived mean first passage times. MatLab is utilized to carry out the computations, using some test problems from the literature.

## 1. Introduction

A variety of techniques have been developed for computing stationary distributions and mean first passage times (MFPTs) of Markov chains (MCs). In this paper we focus primarily on the accurate computation of these key properties based upon the well known state reduction procedures of Grassman, Taksar and Heyman (the GTH algorithm) ([4]), or the equivalent Sheskin ([15]) procedure, that were developed primarily for the computation of the stationary distributions of irreducible MCs. The stability of the procedure is the result of the observation that no subtractions need be carried out. This is discussed in Section 2. Kohlas [13] developed a related procedure for the computation of the MFPTs, based mainly on considering the computation of the mean times to absorption, by showing that the computations were more naturally focused on considering the underlying model as a Markov renewal process (MRP) rather than as a MC. We delve into these procedures in more detail after first summarizing, in Section 3, the key properties of MRPs. In Section 4 we work through the ideas of the Kohlas algorithm and give a general procedure for computing the mean passage times between any two states rather than consider mean times to absorption, as in Kohlas [13]. We explore in some detail, in Sections 5 to 8, procedures for the special cases of particular finite state spaces of one, two, three and four states, obtaining expressions for the MFPTs, some of which that have previously been given in the literature. We see that the Kohlas procedure is not ideal for the global derivation of the MFPT matrix but we develop in Section 9 a modification of the Kohlas procedure, an Extended GTH procedure that will lead to expressions for the MFPTs using effectively the same calculations as in the GTH algorithm. In the final section we explore some calculations, using MatLab, of the key properties for some ill conditioned transition matrices that have previously been considered as test problems in the literature.

Note, that due to space considerations, we do not explore other procedures for finding MFPTs. We leave this for a further paper where we compare the procedure of this paper with the well known approaches of Kemeny and Snell [12] and others, for example, in Meyer [14], Stewart [16], Hunter [9] and Dayar and Akar [3] as well as some new perturbation procedures under development by the author. This also enables us to make additional comparisons using the test problems of this paper to validate the stability of this new procedure.

## 2. Computation of the stationary probabilities

Let $P^{(N)} = \left[ p_{ij}^{(N)} \right] = \begin{bmatrix} Q_{N-1}^{(N)} & \boldsymbol{p}_{N-1}^{(N)(c)} \\ \boldsymbol{p}_{N-1}^{(N)(r)T} & p_{NN}^{(N)} \end{bmatrix}$ be the $N \times N$ transition matrix associated with an irreducible MC $\{X_k^{(N)}, k \geq 0\}$ with state space $S_N = \{1, 2, ..., N\}$ of $N$ states.

Let $\boldsymbol{p}_{N-1}^{(N)(r)T} = (p_{N,1}^{(N)}, p_{N,2}^{(N)} ..., p_{N,N-1}^{(N)})$ and $\boldsymbol{p}_{N-1}^{(N)(c)T} = (p_{1,N}^{(N)}, ..., p_{N-1,N}^{(N)})$ be $1 \times (N-1)$ row vectors with $r$ and $c$ denoting, respectively, row and column elements of the probabilities, with the superscript $N$ denoting that they are from the $N$-th row or $N$-th column of the $P^{(N)}$ matrix and the subscript $N-1$ that they are vectors of length $N-1$. Similarly, we use the superscript $N$ in the sub-matrices $Q_{N-1}^{(N)}$ to denote that they are sub-matrices of the transition matrix $P^{(N)}$ associated with an $N$-state MC and the subscript $N-1$ to denote that the matrix is of order $(N-1) \times (N-1)$.

Let $\boldsymbol{e}^{(N)T} = (1, 1, ..., 1)$ be an $1 \times N$ vector and $I_N$ be the $N \times N$ identity matrix.

In the procedures that we consider for finding stationary distributions and the MFPTs of MCs, we start with an $N$–state MC $\{X_k^{(N)}, k \geq 0\}$ and reduce the state space by one state at a time. Once we get to two states we expand the state space one state at a time until we return to the final set of $N$ states. We concentrate on the sequential state reduction process at first by starting with $N$ states $1, 2, ..., N$ and initially reducing the state space to $1, 2, ..., N-1$.

For simplicity, when there is no ambiguity, we write $p_{ij}^{(N)}$ simply as $p_{ij}$,

Note that $Q_{N-1}^{(N)}$ is not stochastic, since

$$P^{(N)}\boldsymbol{e}^{(N)} = \begin{bmatrix} Q_{N-1}^{(N)} & \boldsymbol{p}_{N-1}^{(N)(c)} \\ \boldsymbol{p}_{N-1}^{(N)(r)T} & p_{NN}^{(N)} \end{bmatrix} \begin{bmatrix} \boldsymbol{e}^{(N-1)} \\ 1 \end{bmatrix} = \begin{bmatrix} Q_{N-1}^{(N)}\boldsymbol{e}^{(N-1)} + \boldsymbol{p}_{N-1}^{(N)(c)} \\ \boldsymbol{p}_{N-1}^{(N)(r)T}\boldsymbol{e}^{(N-1)} + p_{NN}^{(N)} \end{bmatrix} = \begin{bmatrix} \boldsymbol{e}^{(N-1)} \\ 1 \end{bmatrix}$$

implying that 
$$Q_{N-1}^{(N)}\boldsymbol{e}^{(N-1)} + \boldsymbol{p}_{N-1}^{(N)(c)} = \boldsymbol{e}^{(N-1)} \quad (1)$$

and that 
$$\boldsymbol{p}_{N-1}^{(N)(r)T}\boldsymbol{e}^{(N-1)} + p_{NN}^{(N)} = 1. \quad (2)$$

Let $\boldsymbol{\pi}^{(N)T} = (\pi_1^{(N)}, \pi_2^{(N)}, ...., \pi_{N-1}^{(N)}, \pi_N^{(N)})$ be the stationary probability vector of the $N$-state MC $\{X_k^{(N)}, k \geq 0\}$ with transition matrix $P^{(N)}$ so that

$$\boldsymbol{\pi}^{(N)T} = \boldsymbol{\pi}^{(N)T} P^{(N)}. \tag{3}$$

Let $\boldsymbol{\rho}^{(N-1)T} = (\rho_1^{(N-1)}, \rho_2^{(N-1)}, ...., \rho_{N-1}^{(N-1)}) = (\pi_1^{(N)}, \pi_2^{(N)}, ...., \pi_{N-1}^{(N)})$ so that $\boldsymbol{\pi}^{(N)T} = (\boldsymbol{\rho}^{(N-1)T}, \pi_N^{(N)})$.

From (3),

$$\boldsymbol{\pi}^{(N)T} = (\boldsymbol{\rho}^{(N-1)T}, \pi_N^{(N)}) = (\boldsymbol{\rho}^{(N-1)T}, \pi_N^{(N)}) \begin{bmatrix} Q_{N-1}^{(N)} & \boldsymbol{p}_{N-1}^{(N)(c)} \\ \boldsymbol{p}_{N-1}^{(N)(r)T} & p_{NN}^{(N)} \end{bmatrix}$$

$$= (\boldsymbol{\rho}^{(N-1)T} Q_{N-1}^{(N)} + \pi_N^{(N)} \boldsymbol{p}_{N-1}^{(N)(r)T}, \boldsymbol{\rho}^{(N-1)T} \boldsymbol{p}_{N-1}^{(N)(c)} + \pi_N^{(N)} p_{NN}^{(N)}),$$

implying that $\boldsymbol{\rho}^{(N-1)T} Q_{N-1}^{(N)} + \pi_N^{(N)} \boldsymbol{p}_{N-1}^{(N)(r)T} = \boldsymbol{\rho}^{(N-1)T}$ (4)

and $\boldsymbol{\rho}^{(N-1)T} \boldsymbol{p}_{N-1}^{(N)(c)} + \pi_N^{(N)} p_{NN}^{(N)} = \pi_N^{(N)}.$ (5)

Equations (2) and (5) imply that

$$\pi_N^{(N)} = \frac{\boldsymbol{\rho}^{(N-1)T} \boldsymbol{p}_{N-1}^{(N)(c)}}{\boldsymbol{p}_{N-1}^{(N)(r)T} \boldsymbol{e}^{(N)}} = \frac{\sum_{i=1}^{N-1} \rho_i^{(N-1)} p_{iN}^{(N)}}{\sum_{i=1}^{N-1} p_{Nj}^{(N)}}, \tag{6}$$

expressing $\pi_N^{(N)}$ in terms of $\rho_1^{(N-1)}, \rho_2^{(N-1)}, ...., \rho_{N-1}^{(N-1)}$ and the transition probabilities associated with $P^{(N)}$.

Further, from equations (4) and (6),

$$\boldsymbol{\rho}^{(N-1)T} \left( I_{N-1} - Q_{N-1}^{(N)} - \frac{\boldsymbol{p}_{N-1}^{(N)(c)} \boldsymbol{p}_{N-1}^{(N)(r)T}}{\boldsymbol{p}_{N-1}^{(N)(r)T} \boldsymbol{e}^{(N-!)}} \right) = \boldsymbol{0}^T. \tag{7}$$

Let $$P^{(N-1)} = Q_{N-1}^{(N)} + \frac{\boldsymbol{p}_{N-1}^{(N)(c)} \boldsymbol{p}_{N-1}^{(N)(r)T}}{\boldsymbol{p}_{N-1}^{(N)(r)T} \boldsymbol{e}^{(N-1)}}. \tag{8}$$

Note that $P^{(N-1)}$ is a stochastic matrix with $N$–1 states, since from (1),

$$P^{(N-1)} \boldsymbol{e}^{(N-1)} = Q_{N-1}^{(N)} \boldsymbol{e}^{(N-1)} + \frac{\boldsymbol{p}_{N-1}^{(N)(c)} \boldsymbol{p}_{N-1}^{(N)(r)T} \boldsymbol{e}^{(N-1)}}{\boldsymbol{p}_{N-1}^{(N)(r)T} \boldsymbol{e}^{(N-1)}} = \boldsymbol{e}^{(N-1)} - \boldsymbol{p}_{N-1}^{(N)(c)} + \boldsymbol{p}_{N-1}^{(N)(c)} = \boldsymbol{e}^{(N-1)}.$$

Let $\{X_k^{(N-1)}, k \geq 0\}$ be the MC that has $P^{(N-1)}$ as its transition matrix. Note also that $\boldsymbol{p}_{N-1}^{(N)(c)} \boldsymbol{p}_{N-1}^{(N)(r)T}$ is an $(N–1) \times (N–1)$ matrix whose $(i,j)$-th element is $p_{iN}^{(N)} p_{Nj}^{(N)}$, so that if we write $P^{(N)} = \left[ p_{ij}^{(N)} \right]$ with $P^{(N-1)} = \left[ p_{ij}^{(N-1)} \right]$ then, from (8),

$$p_{ij}^{(N-1)} = p_{ij}^{(N)} + \frac{p_{iN}^{(N)} p_{Nj}^{(N)}}{S(N)}, \quad 1 \leq i \leq N-1, \ 1 \leq j \leq N-1; \tag{9}$$

where $S(N) \equiv 1 - p_{NN}^{(N)} = \sum_{j=1}^{N-1} p_{Nj}^{(N)} = \boldsymbol{p}_{N-1}^{(N)(r)T} \boldsymbol{e}^{(N-1)}$ since $P^{(N)}$ is a stochastic matrix.
Note the computation of the quantities $S(N)$ can be carried out without any subtraction.

We can interpret the transition probabilities $p_{ij}^{(N-1)}$ in the MC $\{X_k^{(N-1)}, k \geq 0\}$ on the state space $S_{N-1}$ as the transition probability from state $i$ to $j$ of the MC $\{X_k^{(N)}, k \geq 0\}$ on $S_N$ *restricted to* $S_{N-1}$, i.e. the "*censored*" MC. ( See also pg. 17, Bini, Latouche and Meini [2]).
For $(i,j) \in S_{N-1} \times S_{N-1}$ it is possible to jump directly from $i$ to $j$ with probability $p_{ij}^{(N)}$. Alternatively, it is possible to jump from $i$ to $j$ via state $N$, being held at state $N$ for $t$ steps, ($t$ = 0, 1, 2,…) followed by a one-step jump to $j$ from $N$, with probability

$$p_{iN}^{(N)}\left(\sum_{i=0}^{\infty}\left(p_{NN}^{(N)}\right)^{t}\right)p_{Nj}^{(N)} = \frac{p_{iN}^{(N)} p_{Nj}^{(N)}}{1-p_{NN}^{(N)}} = \frac{p_{iN}^{(N)} p_{Nj}^{(N)}}{S(N)},$$ leading to the general expression (9) for $p_{ij}^{(N-1)}$.

Note that there is a connection between equation (9) and Schur complementation. This is discussed in Bini, Latouche and Meini, (pg 17, [2]).

Note that if the MC $\{X_k^{(N)}, k \geq 0\}$ with state space $S_N$ is irreducible (i.e. each state $j$ can be reached from state $i$ in a finite number of $k$ steps) then the MC $\{X_k^{(N-1)}, k \geq 0\}$ with state space $S_{N-1}$ is also irreducible since there will still be a path from state $j$ that can be reached from state $i$ in either the same $k$ steps, if avoiding state $N$, or in a fewer number of steps if passing through $N$ in the original MC $\{X_k^{(N)}, k \geq 0\}$.

Further, from (7) and (8), $\boldsymbol{\rho}^{(N-1)T}\left(I_{N-1} - P^{(N-1)}\right) = \boldsymbol{0}^T$, so that $\boldsymbol{\rho}^{(N-1)T}$ satisfies the property required for a stationary probability vector of the irreducible MC $\{X_k^{(N-1)}, k \geq 0\}$ on $S_{N-1}$ with transition matrix $P^{(N-1)}$, and hence

$$\boldsymbol{\rho}^{(N-1)T} = c_{N-1} \boldsymbol{\pi}^{(N-1)T}, \tag{10}$$

where $\boldsymbol{\pi}^{(N-1)T} \boldsymbol{e}^{(N-1)} = 1$, implying that $c_{N-1} = \boldsymbol{\rho}^{(N-1)T} \boldsymbol{e}^{(N-1)}$.

Note that $c_{N-1} = \sum_{i=1}^{N-1} \pi_i^{(N)} = 1 - \pi_N^{(N)}$. Thus

$$\boldsymbol{\pi}^{(N-1)T} = (\pi_1^{(N-1)}, \pi_2^{(N-1)}, \ldots, \pi_{N-1}^{(N-1)}) = \frac{\boldsymbol{\rho}^{(N-1)T}}{\boldsymbol{\rho}^{(N-1)T} \boldsymbol{e}^{(N-1)}} = \frac{(\pi_1^{(N)}, \pi_2^{(N)}, \ldots, \pi_{N-1}^{(N)})}{1 - \pi_N^{(N)}},$$

and hence 
$$\pi_i^{(N-1)} = \frac{\pi_i^{(N)}}{1 - \pi_N^{(N)}} = \frac{\pi_i^{(N)}}{\sum_{k=1}^{N-1} \pi_k^{(N)}}, \quad 1 \leq i \leq N-1. \tag{11}$$

Thus we have reduced the state space from $N$ to $N-1$ with the resulting MC $\{X_k^{(N-1)}, k \geq 0\}$ having a stationary distribution $\{\pi_i^{(N-1)}\}$ that is a scaled version of the first $N-1$ components of the stationary distribution of the MC $\{X_k^{(N)}, k \geq 0\}$ with $N$ states, as given by (11).

Let us define the stationary probability vector of the MC $\{X_k^{(N)}, k \geq 0\}$ as $\boldsymbol{\pi}^T = (\pi_1, \pi_2, \ldots, \pi_N) = \boldsymbol{\pi}^{(N)T}$. As we continue to reduce the state space to $S_n$ ($n = 1, 2, \ldots, N-1$) it is clear, from an extension of (11), that

$$\boldsymbol{\pi}^{(n)T} = (\pi_1^{(n)}, \pi_2^{(n)}, \ldots, \pi_n^{(n)}) = k_n(\pi_1, \pi_2, \ldots, \pi_n) \text{ where } k_n = 1 \Big/ \sum_{i=1}^{n} \pi_i. \tag{12}$$

i.e. the stationary probabilities of the MC $\{X_k^{(n)}, k \geq 0\}$ on $S_n$ are scaled versions of the first $n$ stationary probabilities of the MC $\{X_k^{(N)}, k \geq 0\}$ on $S_N$.

Let us now consider expanding the state space from $S_{N-1}$ to $S_N$. Note that, from (11),

$$\pi_i^{(N)} = (1 - \pi_N^{(N)})\pi_i^{(N-1)} = c_{N-1} \pi_i^{(N-1)}, \quad 1 \leq i \leq N-1. \tag{13}$$

i.e. the first $N$-1 terms of $\pi_i^{(N)}$ are a multiple of $\pi_i^{(N-1)}$. Further, from (6) and definition of $S(N)$,

$$\pi_N^{(N)} = c_{N-1} \frac{\sum_{i=1}^{N-1} \pi_i^{(N-1)} p_{iN}^{(N)}}{S(N)}, \tag{14}$$

From (13) and (14), the constant $c_{N-1}$ is determined from the fact that $\sum_{i=1}^{N} \pi_i^{(N)} = 1$, and the stationary distribution for the MC on $S_N$ can be determined from the MC on $S_{N-1}$ yielding

$$\boldsymbol{\pi}^{(N)T} = c_{N-1}(\pi_1^{(N-1)}, \pi_2^{(N-1)}, \ldots, \pi_{N-1}^{(N-1)}, \frac{\sum_{i=1}^{N-1} \pi_i^{(N-1)} p_{iN}^{(N)}}{S(N)}) . \tag{15}$$

leading to a procedure for determining the stationary distribution on the expanded state space.

The reduction process continues until we reach the state space $S_2 = \{1,2\}$ when we obtain the irreducible stochastic matrix $P^{(2)} = \begin{bmatrix} p_{11}^{(2)} & p_{12}^{(2)} \\ p_{21}^{(2)} & p_{22}^{(2)} \end{bmatrix}$ associated with the MC $\{X_k^{(2)}, k \geq 0\}$.

The stationary probability vector of this MC is given by $\boldsymbol{\pi}^{(2)T} = (\pi_1^{(2)}, \pi_2^{(2)})$.

The second stationary equation is $\pi_2^{(2)} = \pi_1^{(2)} p_{12}^{(2)} + \pi_2^{(2)} p_{22}^{(2)}$ implying $\pi_2^{(2)}(1 - p_{22}^{(2)}) = \pi_1^{(2)} p_{12}^{(2)}$, i.e. 
$$\pi_2^{(2)} = \pi_1^{(2)} \frac{p_{12}^{(2)}}{1 - p_{22}^{(2)}} = \pi_1^{(2)} \frac{p_{12}^{(2)}}{p_{21}^{(2)}} = \pi_1^{(2)} \frac{p_{12}^{(2)}}{S(2)} . \tag{16}$$

Note that $S(2) = 1 - p_{22}^{(2)} = \sum_{j=1}^{1} p_{2j}^{(2)} = p_{21}^{(2)} = \boldsymbol{p}_1^{(2)(r)T} \boldsymbol{e}^{(1)}$.

Since from (12), $\boldsymbol{\pi}^{(2)T} = (\pi_1^{(2)}, \pi_2^{(2)}) = k_2(\pi_1, \pi_2)$, we have from (16) by dividing by $k_2$ that

$$\pi_2 = \pi_1 \frac{p_{12}^{(2)}}{S(2)}. \tag{17}$$

We now proceed with increasing the state space using the process described above.
Observe that from (12) with $n = 3$, and (15) with $N = 3$,

$$\boldsymbol{\pi}^{(3)T} = (\pi_1^{(3)}, \pi_2^{(3)}, \pi_3^{(3)}) = k_3(\pi_1, \pi_2, \pi_3) = c_2(\pi_1^{(2)}, \pi_2^{(2)}, \frac{\sum_{i=1}^{2} \pi_i^{(2)} p_{i3}^{(3)}}{S(3)})$$

implying $\pi_3^{(3)} = \frac{\sum_{i=1}^{2} \pi_i^{(2)} p_{i3}^{(3)}}{\sum_{i=1}^{2} p_{3i}^{(3)}} = \frac{\pi_1^{(2)} p_{13}^{(3)} + \pi_2^{(2)} p_{23}^{(3)}}{1 - p_{33}^{(3)}} = \pi_1^{(2)} \frac{p_{13}^{(3)}}{S(3)} + \pi_2^{(2)} \frac{p_{23}^{(3)}}{S(3)}$, and hence by scaling

$$\pi_3 = \pi_1 \frac{p_{13}^{(3)}}{S(3)} + \pi_2 \frac{p_{23}^{(3)}}{S(3)}. \tag{18}$$

leading in general, for $n = 2, \ldots, N$, to

$$\pi_n = \frac{\sum_{i=1}^{n-1} \pi_i p_{in}^{(n)}}{\sum_{i=1}^{n-1} p_{ni}^{(n)}} = \sum_{i=1}^{n-1} \pi_i \frac{p_{in}^{(n)}}{S(n)}. \tag{19}$$

Thus if $\pi_j = kr_j$ with $r_1 = 1$ then $\sum_{i=1}^{N} \pi_i = 1 \Rightarrow k = 1 / \sum_{i=1}^{N} r_i$ with

$$r_n = \frac{\sum_{i=1}^{n-1} r_i p_{in}^{(n)}}{S(n)}, (n = 2, \ldots, N), \text{ implying } \pi_i = r_i / \sum_{n=1}^{N} r_n, i = 1, 2, \ldots, N.$$

We summarize the procedure as follows.

**Theorem 1:**

Given a finite irreducible MC $\{X_k^{(N)}, k \geq 0\}$ with state space $S_N = \{1, 2, \ldots, N\}$ and transition matrix $P^{(N)} = \left[p_{ij}^{(N)}\right]$ its stationary probabilities $\{\pi_i^{(N)}\}$ can be computed as follows:

1. Compute, successively for $n = N, N-1, \ldots, 3$,

$$p_{ij}^{(n-1)} = p_{ij}^{(n)} + \frac{p_{in}^{(n)} p_{nj}^{(n)}}{S(n)}, \quad 1 \leq i \leq n-1, \ 1 \leq j \leq n-1; \text{ where } S(n) = \sum_{j=1}^{n-1} p_{nj}^{(n)}.$$

2. Set $r_1 = 1$ and compute successively for $n = 2, \ldots, N$, $r_n = \dfrac{\sum_{i=1}^{n-1} r_i p_{in}^{(n)}}{S(n)}$.

3. Compute, for $i = 1, 2, \ldots, N$, $\pi_i^{(N)} = \dfrac{r_i}{\sum_{j=1}^{N} r_j}$.

This is a formal derivation of the Grassman, Taksar and Heyman (GTH) algorithm ([4]) or the equivalent Sheskin State Reduction procedure ([15]) for finding the stationary distribution of an irreducible finite MC. The procedure is numerically stable and accurate in that no subtractions need be carried out.

Note that, as a result of equation (12), the stationary distribution for the derived MC $\{X_k^{(n)}, k \geq 0\}$ with transition probability matrix $P^{(n)} = \left[p_{ij}^{(n)}\right]$ on the reduced state space $S_n$, is given as $\pi_i^{(n)} = \dfrac{r_i}{\sum_{j=1}^{n} r_j}$, $i = 1, 2, \ldots, n$.

## 3. Markov renewal processes

We give a brief review of some of the key properties of Markov renewal processes. We refer the reader to Section 2.2 of [9] where the following general concepts and notation are presented.

We consider a MRP $\{(X_n, T_n), n \geq 0\}$, with state space $S = \{1, 2, \ldots, N\}$ and semi-Markov kernel $Q(t) = \left[Q_{ij}(t)\right]$, where $Q_{ij}(t) = P\{X_{n+1} = j, T_{n+1} - T_n \leq t | X_n = i\}$, $(i, j) \in S$. $\{X_n\}$, $(n \geq 0)$, tracks the states successively visited and $T_n$ is the time of the *n*-th transition. Observe that $Q_{ij}(+\infty) = P\{X_{n+1} = j | X_n = i\}$, so that $\{X_n\}$ is a MC, the embedded MC, with transition matrix $P = \left[p_{ij}\right]$ where $p_{ij} = Q_{ij}(+\infty)$. Further, we can express $Q_{ij}(t)$ as $Q_{ij}(t) = p_{ij} F_{ij}(t)$ where $F_{ij}(t) = P\{T_{n+1} - T_n \leq t | X_n = i, X_{n+1} = j\}$. Thus $F_{ij}(t)$ is the distribution function of the "*holding time*" $T_{n+1} - T_n$ in state $X_n$ until transition into state $X_{n+1}$ given that the MRP makes a transition from $X_n$ to $X_{n+1}$.

Let $\mu_{ij} = \int_0^\infty t \, dQ_{ij}(t)$ so that $\mu_{ij} = p_{ij} E[T_{n+1} - T_n | X_n = i, X_{n+1} = j]$.

We assume that the embedded MC $\{X_n, n \geq 0\}$ is irreducible and hence has a stationary distribution $\{\pi_j\}, (j \in S)$ and associated stationary probability vector $\pi^T = (\pi_1, \pi_2, ...., \pi_N)$. Let $N = [\mu_{ij}]$ then, (equation (2.10), [9]), the MFPT matrix $M = [m_{ij}]$ of the MRP $\{(X_n, T_n), n \geq 0\}$ satisfies the equation

$$(I - P)M = NE - PM_d, \qquad (20)$$

where $M_d = [\delta_{ij} m_{ij}] = diag(m_{11}, ...., m_{NN})$, (with $\delta_{ij} = 1$ when $i = j$ and 0 otherwise.)

Let $\mu = Ne$ so that $NE = Nee^T = \mu e^T$. If $\mu^T = (\mu_1, \mu_2, ..., \mu_N)$ then $\mu_i = \sum_{j=1}^{N} \mu_{ij}$. Observe that $\mu_i = E[T_{n+1} - T_n | X_n = i]$, the "*expected holding time starting in state i*". We introduce one further piece of notation. The "*mean asymptotic increment*", for the MRP is given by $\lambda_1 = \pi^T \mu$, i.e. the "*expected holding time under stationary conditions*". From Section 5.2 of [9], $M_d = \lambda_1 (\Pi_d)^{-1}$ where $\Pi = e\pi^T$ implying that

$$m_{jj} = \frac{\lambda_1}{\pi_j}. \qquad (21)$$

Note that when $T_{n+1} = T_n + 1$, the MRP $\{(X_n, T_n), n \geq 0\}$ reduces to a discrete time MC $\{X_n, n \geq 0\}$ with $\mu_{ij} = p_{ij}$, $\mu_i = 1$, for all $i$ and $\lambda_1 = 1$. Thus $\mu = e$ and $NE = ee^T = E$.

## 4. Computation of the Mean First Passage Times

We seek a computational procedure that will enable us to calculate all the MFPTs times in a MC.

As Kohlas [13] pointed out in his pioneering paper, it is more natural to consider the Markov renewal setting. Let us define $M_n = [m_{ij}]$, $(1 \leq i \leq n, 1 \leq j \leq n)$ as the MFPT matrix of the MRP $\{(X_k^{(n)} T_k^{(n)}), k \geq 0\}$ with $n$-states, transition matrix $P^{(n)}$ and mean holding time vector $\mu^{(n)}$. From (20), the matrix $M_n$ satisfies

$$(I_n - P^{(n)})M_n = \mu^{(n)} e^{(n)T} - P^{(n)}(M_n)_d. \qquad (22)$$

Note that for the MC $\{X_k^{(N)}, k \geq 0\}$ starting with $N$ states, $\mu^{(N)T} = e^{(N)T} = (1,1,...,1)$.

Let us partition $M_n$ as $M_n = \begin{bmatrix} M_{n-1} & m_{n-1}^{(n)(c)} \\ m_{n-1}^{(n)(r)T} & m_{nn} \end{bmatrix} \qquad (23)$

where $M_{n-1} = [m_{ij}]$, $(1 \leq i \leq n-1, 1 \leq j \leq n-1)$, $m_{n-1}^{(n)(r)T} = (m_{n1}, m_{n2}, ...., m_{n,n-1})$ and $m_{n-1}^{(n)(c)T} = (m_{1n}, m_{2n}, ...., m_{n-1,n})$.

Define $\mu^{(n)T} = (\mu_1^{(n)}, ..., \mu_{n-1}^{(n)}, \mu_n^{(n)}) = (\mu_{n-1}^{(n)T}, \mu_n^{(n)})$ where $\mu_{n-1}^{(n)T} = (\mu_1^{(n)}, ..., \mu_{n-1}^{(n)})$. We partition $P^{(n)} = \begin{bmatrix} Q_{n-1}^{(n)} & p_{n-1}^{(n)(c)} \\ p_{n-1}^{(n)(r)T} & p_{nn}^{(n)} \end{bmatrix}$ so that block multiplication of (22) yields

$$\begin{bmatrix} I_{n-1} - Q_{n-1}^{(n)} & -p_{n-1}^{(n)(c)} \\ -p_{n-1}^{(n)(r)T} & 1 - p_{nn}^{(n)} \end{bmatrix} \begin{bmatrix} M_{n-1} & m_{n-1}^{(n)(c)} \\ m_{n-1}^{(n)(r)T} & m_{nn} \end{bmatrix}$$

$$= \begin{bmatrix} \boldsymbol{\mu}_{n-1}^{(n)} e^{(n-1)T} & \boldsymbol{\mu}_{n-1}^{(n)} \\ \mu_n^{(n)} e^{(n-1)T} & \mu_n^{(n)} \end{bmatrix} - \begin{bmatrix} Q_{n-1}^{(n)} & \boldsymbol{p}_{n-1}^{(n)(c)} \\ \boldsymbol{p}_{n-1}^{(n)(r)T} & p_{nn}^{(n)} \end{bmatrix} \begin{bmatrix} (M_{n-1})_d & \boldsymbol{0} \\ \boldsymbol{0}^T & m_{nn} \end{bmatrix}.$$

Hence

(1,1) Block: $(I_{n-1} - Q_{n-1}^{(n)})M_{n-1} - \boldsymbol{p}_{n-1}^{(n)(c)} \boldsymbol{m}_{n-1}^{(n)(r)T} = \boldsymbol{\mu}_{n-1}^{(n)} e^{(n-1)T} - Q_{n-1}^{(n)}(M_{n-1})_d$. (24)

(1,2) Block: $(I_{n-1} - Q_{n-1}^{(n)})\boldsymbol{m}_{n-1}^{(n)(c)} - m_{nn}\boldsymbol{p}_{n-1}^{(n)(c)} = \boldsymbol{\mu}_{n-1}^{(n)} - m_{nn}\boldsymbol{p}_{n-1}^{(n)(c)}$. (25)

(2,1) Block: $-\boldsymbol{p}_{n-1}^{(n)(r)T} M_{n-1} + (1 - p_{nn}^{(n)})\boldsymbol{m}_{n-1}^{(r)T} = \mu_n^{(n)} e^{(n-1)T} - \boldsymbol{p}_{n-1}^{(n)(r)T}(M_{n-1})_d$. (26)

(2,2) Block: $-\boldsymbol{p}_{n-1}^{(n)(r)T} \boldsymbol{m}_{n-1}^{(n)(c)} + (1 - p_{nn}^{(n)})m_{nn} = \mu_n^{(n)} - p_{nn}^{(n)} m_{nn}$. (27)

From (26),

$$\boldsymbol{m}_{n-1}^{(n)(r)T} = \frac{1}{(1 - p_{nn}^{(n)})}\left\{\boldsymbol{p}_{n-1}^{(n)(r)T}\left(M_{n-1} - (M_{n-1})_d\right) + \mu_n^{(n)} e^{(n-1)T}\right\}$$

and, using (2),

$$\boldsymbol{m}_{n-1}^{(n)(r)T} = \frac{1}{\boldsymbol{p}_{n-1}^{(n)(r)T} e^{(n-1)}}\left\{\boldsymbol{p}_{n-1}^{(n)(r)T}\left(M_{n-1} - (M_{n-1})_d\right) + \mu_n^{(n)} e^{(n-1)T}\right\}. \quad (28)$$

Substitute into (24)

$$\left(I_{n-1} - Q_{n-1}^{(n)} - \frac{\boldsymbol{p}_{n-1}^{(n)(c)} \boldsymbol{p}_{n-1}^{(n)(r)T}}{\boldsymbol{p}_{n-1}^{(n)(r)T} e^{(n-1)}}\right)M_{n-1} = \left(\boldsymbol{\mu}_{n-1}^{(n)} + \frac{\mu_n^{(n)} \boldsymbol{p}_{n-1}^{(n)(c)}}{\boldsymbol{p}_{n-1}^{(n)(r)T} e^{(n-1)}}\right)e^{(n-1)T} + \left(Q_{n-1}^{(n)} + \frac{\boldsymbol{p}_{n-1}^{(n)(c)} \boldsymbol{p}_{n-1}^{(n)(r)T}}{\boldsymbol{p}_{n-1}^{(n)(r)T} e^{(n-1)}}\right)(M_{n-1})_d.$$

Thus, using the expression for $P^{(n-1)}$ as derived earlier (cf. equation (8)),

$$(I_{n-1} - P^{(n-1)})M_{n-1} = \boldsymbol{\mu}^{(n-1)} e^{(n-1)T} - P^{(n-1)}(M_{n-1})_d \quad \text{where} \quad \boldsymbol{\mu}^{(n-1)} = \boldsymbol{\mu}_{n-1}^{(n)} + \frac{\mu_n^{(n)} \boldsymbol{p}_{n-1}^{(n)(c)}}{\boldsymbol{p}_{n-1}^{(n)(r)T} e^{(n-1)}}. \quad (29)$$

This is of similar form to the $n$-state case as given by (22) but with the state space reduced to $n - 1$ and a changed form for $\boldsymbol{\mu}^{(n-1)}$.

This leads to the following structural result.

**Theorem 2:**
Let $\{(X_k^{(n)}, T_k^{(n)}), k \geq 0\}$ be a MRP with state space $S_n = \{1, 2, \ldots, n\}$, ($n = 2, \ldots, N$), transition matrix $P^{(n)} = [p_{ij}^{(n)}]$, MFPT matrix $M_n = [m_{ij}]$, ($1 \leq i \leq n, 1 \leq j \leq n$), and vector of mean holding times $\boldsymbol{\mu}^{(n)T} = (\mu_1^{(n)}, \ldots, \mu_{n-1}^{(n)}, \mu_n^{(n)})$ then $M_n$ satisfies equation (22), i.e. $(I_n - P^{(n)})M_n = \boldsymbol{\mu}^{(n)} e^{(n)T} - P^{(n)}(M_n)_d$, or, in element form, $m_{ij} = \mu_i^{(n)} + \sum_{k=1, k \neq j}^{n} p_{ik}^{(n)} m_{kj}$, $(1 \leq i \leq n, 1 \leq j \leq n)$. (30)

Then, under the state reduction process as carried out under the GTH algorithm, $\{(X_k^{(n-1)}, T_k^{(n-1)}), k \geq 0\}$ is also a MRP with state space $S_{n-1} = \{1, 2, \ldots, n-1\}$, transition matrix $P^{(n-1)} = [p_{ij}^{(n-1)}]$, and MFPT matrix $M_{n-1} = [m_{ij}]$, $(1 \leq i \leq n-1, 1 \leq j \leq n-1)$, which satisfies equation (29) i.e. $(I_{n-1} - P^{(n-1)})M_{n-1} = \boldsymbol{\mu}^{(n-1)} e^{(n-1)T} - P^{(n-1)}(M_{n-1})_d$,

where the transition probabilities $p_{ij}^{(n-1)}$ are given by

$$p_{ij}^{(n-1)} = p_{ij}^{(n)} + \frac{p_{in}^{(n)} p_{nj}^{(n)}}{S(n)}, \quad 1 \leq i \leq n-1, \, 1 \leq j \leq n-1, \tag{31}$$

and the elements of the mean holding time vector $\boldsymbol{\mu}^{(n-1)T} = (\mu_1^{(n-1)},...,\mu_{n-1}^{(n-1)})$ are given by

$$\mu_i^{(n-1)} = \mu_i^{(n)} + \frac{p_{in}^{(n)} \mu_n^{(n)}}{S(n)}, \quad 1 \leq i \leq n-1, \tag{32}$$

where $S(n) = \boldsymbol{p}_{n-1}^{(n)(r)T} \boldsymbol{e}^{(n-1)} = \sum_{j=1}^{n-1} p_{nj}^{(n)} = 1 - p_{nn}^{(n)}$.

□

Note that equation (31) is identical to format of the transition probabilities as used in the GTH algorithm with the derivation given by equations (8) and (9) with $N$ replaced by $n$. The expression for the elemental expressions for the mean holding times (32) follows from (29).

Thus for the reduced MRP $\{(X_k^{(n-1)} T_k^{(n-1)}), k \geq 0\}$ the MFPTs $m_{ij}, (1 \leq i \leq n-1, 1 \leq j \leq n-1)$ are identical to those of the same pairs of states as in the original MRP $\{(X_k^{(n)} T_k^{(n)}), k \geq 0\}$. This means that we can reduce the state space by successive steps retaining the same MFPTs for the reduced state space in the upper block of $M_n$ although the mean holding times in the states are modified, as given by equation (32).

Upon increasing the state space from $S_{n-1}$ to $S_n$, as in the GTH algorithm, we wish to find an expression for the elements of $M_n$, given $M_{n-1} = [m_{ij}], (1 \leq i \leq n-1, 1 \leq j \leq n-1)$. Thus, from (23), we need to find expressions for $\boldsymbol{m}_{n-1}^{(n)(c)}, \boldsymbol{m}_{n-1}^{(n)(r)T}$ and $m_{nn}$.

From the properties of MCs and MRPs the following results for the mean recurrence times $m_{nn}$, can be deduced:

**Theorem 3:**

(1) $m_{nn} = \dfrac{\lambda_1^{(N)}}{\pi_n^{(N)}}$ where $\lambda_1^{(N)} = \boldsymbol{\pi}^{(N)T} \boldsymbol{\mu}^{(N)} = \sum_{k=1}^{N} \pi_k^{(N)} \mu_k^{(N)}$. (33)

(2) $m_{nn} = \dfrac{\lambda_1^{(n)}}{\pi_n^{(n)}}$ where $\lambda_1^{(n)} = \boldsymbol{\pi}^{(n)T} \boldsymbol{\mu}^{(n)} = \sum_{k=1}^{n} \pi_k^{(n)} \mu_k^{(n)}$. (34)

(3) In the MC setting for $\{X_k^{(N)}, k \geq 0\}$, $m_{nn} = 1/\pi_n^{(N)}$. (35)

**Proof**:
Starting with a MRP $\{(X_k^{(N)}, T_k^{(N)}), k \geq 0\}$ on the state space $S_N = \{1,2,...,N\}$, equation (21) implies that $m_{ii} = \dfrac{\lambda_1^{(N)}}{\pi_i^{(N)}}$, where $\lambda_1^{(N)} = \boldsymbol{\pi}^{(N)T} \boldsymbol{\mu}^{(N)}$, leading to (33) for $i = n$ and also to equation (34) since $\{(X_k^{(n)}, T_k^{(n)}), k \geq 0\}$ is also a MRP with mean increment $\lambda_1^{(n)}$. Equation (35) follows since, in the MC setting, $\mu_i^{(N)} = 1$ and $\lambda_1^{(N)} = 1$.

□

**Theorem 4:**

$$m_{nn} = \mu_n^{(n)} + \sum_{k=1}^{n-1} p_{nk}^{(n)} m_{kn}, \quad n = 2,...,N \tag{36}$$

where $m_{11} = \mu_1^{(1)}$.

**Proof:**
Equation (36) follows from an elemental expression of equation (27). The result for $n = 1$ follows from equation (34) as $\pi_1^{(1)} = 1$ and hence $\lambda_1^{(1)} = \mu_1^{(1)}$. □

Theorem 4 gives an additional useful computational procedure for $m_{nn}$. While it does require knowledge of the $m_{in}$ for $i = 1, 2, \ldots, n - 1$, it avoids the calculation of the stationary distribution which is an advantage in the Markov renewal setting. The computation of the $m_{in}$ for $i < n$ requires some additional computational effort as we shall see shortly.

With knowledge of the elements of $M_{n-1}$ expressions for the elements of $\boldsymbol{m}_{n-1}^{(n)(r)T} = (m_{n1}, m_{n2}, \ldots, m_{n,n-1})$ can easily be deduced directly from equation (28).

**Theorem 5:**
$$m_{nj} = \frac{\mu_n^{(n)} + \sum_{k=1, k \neq j}^{n-1} p_{nk}^{(n)} m_{kj}}{S(n)}, \quad j = 1, \ldots, n-1, \tag{37}$$

where $S(n) = 1 - p_{nn}^{(n)} = \sum_{j=1}^{n-1} p_{nj}^{(n)}$. □

Application of Theorem 5 requires retention of the elements $p_{nk}^{(n)}$ of the $n$-th row of $P^{(n)}$.

It is a little more difficult to find the vector $\boldsymbol{m}_{n-1}^{(n)(c)} = (m_{1n}, m_{2n}, \ldots, m_{n-1,n})$.
From (25), $\qquad (I_{n-1} - Q_{n-1}^{(n)}) \boldsymbol{m}_{n-1}^{(n)(c)} = \boldsymbol{\mu}_{n-1}^{(n)}.$ (38)
Even though $(I_{n-1} - Q_{n-1}^{(n)})^{-1}$ exists we use the reduction procedure used above by eliminating $m_{n-1,n}$ from $\boldsymbol{m}_{n-1}^{(n)(c)T}$ and replacing it in the expressions for the elements $m_{1n}, m_{2n}, \ldots, m_{n-2,n}$. The following theorem enables us to develop expressions for the $m_{in}$ for $i < n$.

**Theorem 6:**
(a) $\quad m_{in} = \mu_i^{(n)} + \sum_{k=1}^{n-1} p_{ik}^{(n)} m_{kn}, \quad i = 1, \ldots, n-1, \ n = 2, \ldots, N.$ (39)

(b) $\quad m_{in} = v_i^{(t,n)} + \sum_{k=1}^{t} q_{ik}^{(t,n)} m_{kn}, 1 \leq i \leq t \leq n-1, \ n = 2, \ldots, N.$ (40)

where $q_{ik}^{(t-1,n)} = q_{ik}^{(t,n)} + \dfrac{q_{it}^{(t,n)} q_{tk}^{(t,n)}}{1 - q_{tt}^{(t,n)}}, \quad i,k = 1, \ldots, t-1, \ t = 2, \ldots, n-1;$ (41)

with $q_{ik}^{(n-1,n)} = p_{ik}^{(n)}, \ i,k = 1, \ldots, n-1, \ n = 2, \ldots, N,$ (42)

and $v_i^{(t-1,n)} = v_i^{(t,n)} + \dfrac{q_{it}^{(t,n)} v_t^{(t,n)}}{1 - q_{tt}^{(t,n)}}, \ i = 1, \ldots, t-1, \ t = 2, \ldots, n-1,$ (43)

with $v_i^{(n-1,n)} = \mu_i^{(n)}, \ i = 1, \ldots, n-1; \ n = 2, \ldots, N$. (44)

(c) $m_{1n} = \dfrac{v_1^{(1,n)}}{R(1,n)}, \ n = 2, \ldots, N,$ (45)

where $R(i,n) = 1 - q_{ii}^{(i,n)}, \ i = 1, \ldots, n-1; \ n = 2, \ldots, N.$ (46)

(d) $m_{in} = \dfrac{v_i^{(i,n)} + \sum_{k=1}^{i-1} q_{ik}^{(i,n)} m_{kn}}{R(i,n)}, \ i = 2, \ldots, n-1; \ n = 2, \ldots, N.$ (47)

**Proof**:

(a) Expression (39) is equation (38) in element form, using equations (42) and (44).

(b) In the first instance when $t = n - 1$, expression (40) is identical to (39).
Now from equation (39) express $m_{n-1,n}$ in terms of the $m_{1n}, m_{2n}, ..., m_{n-2,n}$ obtaining

$$m_{n-1,n} = \frac{\mu_{n-1}^{(n)} + \sum_{k=1}^{n-2} p_{n-1,k}^{(n)} m_{kn}}{1 - p_{n-1,n-1}^{(n)}}. \tag{48}$$

Substitute expression (47) for $m_{n-1,n}$ in each of the $m_{in}, (i = 1, ..., n-2)$, expressions given by (39) to obtain, using equations (44) and (45),

$$m_{in} = \left\{ \mu_i^{(n)} + \frac{p_{i,n-1}^{(n)} \mu_{n-1}^{(n)}}{1 - p_{n-1,n-1}^{(n)}} \right\} + \sum_{k=1}^{n-2} \left\{ p_{ik}^{(n)} + \frac{p_{i,n-1}^{(n)} p_{n-1,k}^{(n)}}{1 - p_{n-1,n-1}^{(n)}} \right\} m_{kn}$$

$$= \left\{ v_i^{(n-1,n)} + \frac{q_{i,n-1}^{(n-1,n)} v_{n-1}^{(n-1,n)}}{1 - q_{n-1,n-1}^{(n-1,n)}} \right\} + \sum_{k=1}^{n-2} \left\{ q_{ik}^{(n-1,n)} + \frac{q_{i,n-1}^{(n-1,n)} q_{n-1,k}^{(n-1,n)}}{1 - q_{n-1,n-1}^{(n-1,n)}} \right\} m_{kn}$$

$$= v_i^{(n-2,n)} + \sum_{k=1}^{n-2} q_{ik}^{(n-2,n)} m_{kn}, \quad 1 \leq i \leq n-2,$$

establishing that equation (40) is true for $t = n - 2$.
We now use a proof by induction. Assume that equation (40) is true for $t = s \leq n - 1$.

Thus $m_{sn} = v_s^{(s,n)} + \sum_{k=1}^{s} q_{sk}^{(s,n)} m_{kn}$, implying that $m_{sn} = \dfrac{v_s^{(s,n)} + \sum_{k=1}^{s-1} q_{sk}^{(s,n)} m_{kn}}{1 - q_{ss}^{(s,n)}}$.

Substitution in equation (40) when $t = s$, yields, using equations (41) and (43), that

$$m_{in} = \left\{ v_i^{(s,n)} + \frac{q_{is}^{(s,n)} v_s^{(s,n)}}{1 - q_{ss}^{(s,n)}} \right\} + \sum_{k=1}^{s-1} \left\{ q_{ik}^{(s,n)} + \frac{q_{is}^{(s,n)} q_{sk}^{(s,n)}}{1 - q_{ss}^{(s,n)}} \right\} m_{kn} = v_i^{(s-1,n)} + \sum_{k=1}^{s-1} q_{ik}^{(s-1,n)} m_{kn}.$$

This implies that equation (40) is true for $t = s - 1$. Since equation (40) is true for $t = n - 1$, (by equation (39)) and hence by induction it is true for $t = n - 2, n - 3, ..., 2, 1$.

(c) From equation (40) when $i = t = 1$, $m_{1n} = v_1^{(1,n)} + q_{1k}^{(1,n)} m_{1n}$ leading to equation (45) with the notation of equation (46).

(d) From equation (40), when $i = t = 2$, $m_{2n} = v_2^{(2,n)} + q_{21}^{(2,n)} m_{1n} + q_{22}^{(2,n)} m_{2n}$

so that $m_{2n} = \dfrac{v_2^{(s,n)} + q_{21}^{(2,n)} m_{1n}}{1 - q_{22}^{(2,n)}}$, leading to equation (47) when $i = 2$.

In general, for $i = 2, ..., n - 1$ equation (47) follows directly from equation (40) when $t = i$.
□

Equations (45) and (47) enable successive derivation of $m_{1n}, m_{2n}, ..., m_{n-2,n}, m_{n-1,n}$ following repeated recursion of equation (47) with $i = 1, 2, ..., n - 1$.

In the calculations expressed by equation (47) it would be advantageous if we could express $R(i,n)$ as a sum of terms, with no subtraction, as was the case for the $S(n)$.

Note when $i = n - 1$, equation (47) is equivalent to equation (39) since $q_{n-1,k}^{(n-1,n)} = p_{n-1,k}^{(n)}$ and

$v_i^{(n-1,n)} = \mu_i^{(n)}$ yielding $m_{n-1,n} = \dfrac{v_{n-1}^{(n-1,n)} + \sum_{k=1}^{n-2} q_{n-1,k}^{(n-1,n)} m_{kn}}{R(n-1,n)} = \dfrac{\mu_{n-1}^{(n)} + \sum_{k=1}^{n-2} p_{n-1,k}^{(n)} m_{kn}}{1 - q_{n-1,n-1}^{(n-1,n)}}$,

where $R(n-1,n) = 1 - p_{n-1,n-1}^{(n)} = \sum_{j=1, j \neq n-1}^{n} p_{n-1,j}^{(n)}$, (since $P^{(n)}$ is stochastic), a sum of terms.

When $i = n-2$, $m_{n-2,n} = \dfrac{v_{n-2}^{(n-2,n)} + \sum_{k=1}^{n-3} q_{n-2,k}^{(n-2,n)} m_{kn}}{R(n-2,n)}$,

where $q_{n-2,k}^{(n-2,n)} = q_{n-2,k}^{(n-1,n)} + \dfrac{q_{n-2,n-1}^{(n-1,n)} q_{n-1,k}^{(n-1,n)}}{1 - q_{n-2,n-2}^{(n-1,n)}} = p_{n-2,k}^{(n)} + \dfrac{p_{n-2,n-1}^{(n)} p_{n-1,k}^{(n)}}{1 - p_{n-1,n-1}^{(n)}}$

and $R(n-2,n) \equiv 1 - q_{n-2,n-2}^{(n-2,n)} = 1 - p_{n-2,n-2}^{(n)} - \dfrac{p_{n-2,n-1}^{(n)} p_{n-1,n-2}^{(n)}}{1 - p_{n-1,n-1}^{(n)}}$

$= \dfrac{(1 - p_{n-1,n-1}^{(n)})(1 - p_{n-2,n-2}^{(n)}) - p_{n-2,n-1}^{(n)} p_{n-1,n-2}^{(n)}}{R(n-1,n)}$.

It follows that the numerator of the expression for $R(n–2,n)$ can also be expressed in terms not involving any subtraction since $1 - p_{n-1,n-1}^{(n)} = \sum_{j=1, j \neq n-1, n-2}^{n} p_{n-1,j}^{(n)} + p_{n-1,n-2}^{(n)}$ and $1 - p_{n-2,n-2}^{(n)} = \sum_{j=1, j \neq n-1, n-2}^{n} p_{n-2,j}^{(n)} + p_{n-2,n-1}^{(n)}$.

It is expected that it can be shown that the denominators of the expressions given by (41), (43), (45) and (46), i.e. $R(t,n) = 1 - q_{tt}^{(t,n)}$, can all be expressed in terms not involving subtractions, as we were able to show for the $S(n)$.

The state reduction process can continue to a single state, $n = 1$, where from (30), $m_{11} = \mu_1^{(1)}$. (see Section 5 for a further discussion on this result.)

We can however finish the state reduction process when we are left with $n = 2$ states. From (30), we have four equations

$$m_{11} = \mu_1^{(2)} + p_{12}^{(2)} m_{21}, \quad m_{12} = \mu_1^{(2)} + p_{11}^{(2)} m_{12},$$
$$m_{21} = \mu_2^{(2)} + p_{22}^{(2)} m_{21}, \quad m_{22} = \mu_2^{(2)} + p_{21}^{(2)} m_{12},$$

that are easily solved to yield, using the observation that $1 - p_{11}^{(2)} = p_{12}^{(2)}$ and $1 - p_{22}^{(2)} = p_{21}^{(2)}$,

$$M_2 = \begin{bmatrix} m_{11} & m_{12} \\ m_{21} & m_{22} \end{bmatrix} = \begin{bmatrix} \mu_1^{(2)} + \dfrac{p_{12}^{(2)}}{p_{21}^{(2)}} \mu_2^{(2)} & \dfrac{\mu_1^{(2)}}{p_{12}^{(2)}} \\ \dfrac{\mu_2^{(2)}}{p_{21}^{(2)}} & \dfrac{p_{21}^{(2)}}{p_{12}^{(2)}} \mu_1^{(2)} + \mu_2^{(2)} \end{bmatrix}. \quad (49)$$

Note, from equation (32) with $n = 2$, $m_{11} = \mu_1^{(1)} = \mu_1^{(2)} + \dfrac{p_{12}^{(2)} \mu_2^{(2)}}{S(2)}$, where $S(2) = 1 - p_{22}^{(2)} = p_{21}^{(2)}$, leading to the expression for $m_{11}$ in (49).
Following the state reduction process to $S_2$ we now need to increase the state space to $S_N$ through the inclusion of successive additional states.

From the process outlined in Theorem 2,

$$M_3 = \begin{bmatrix} & M_2 & & m_{13} \\ & & & m_{23} \\ m_{31} & m_{32} & m_{33} \end{bmatrix}, \text{ where the } M_2 \text{ matrix is given by (49).}$$

From Theorem 5, equation (37), $m_{31} = \dfrac{\mu_3^{(3)} + p_{32}^{(3)} m_{21}}{S(3)}$, $m_{32} = \dfrac{\mu_3^{(3)} + p_{31}^{(3)} m_{12}}{S(3)}$.

From Theorem 6, equation (45), $m_{13} = \dfrac{v_1^{(1,3)}}{R(1,3)}$.

From Theorem 6, equation (47), $m_{23} = \dfrac{v_2^{(2,3)} + q_{21}^{(2,3)} m_{13}}{R(2,3)}$.

From Theorem 4, equation (36), $m_{33} = \mu_3^{(3)} + p_{31}^{(3)} m_{13} + p_{32}^{(3)} m_{23}$.

(Alternatively, from equation (34), $m_{33} = \dfrac{\sum_{k=1}^{3} \pi_k^{(3)} \mu_k^{(3)}}{\pi_3^{(3)}}$, but this expression requires derivation of the $\pi_i^{(3)}$ from the first $r_1, r_2$ and $r_3$ terms of the GTH procedure).

$$\text{Thus } M_3 = \begin{bmatrix} m_{11} & m_{12} & \dfrac{v_1^{(1,3)}}{R(1,3)} \\ m_{21} & m_{22} & \dfrac{v_2^{(2,3)} + p_{21}^{(3)} m_{13}}{R(2,3)} \\ \dfrac{\mu_3^{(3)} + p_{32}^{(3)} m_{21}}{S(3)} & \dfrac{\mu_3^{(3)} + p_{31}^{(3)} m_{12}}{S(3)} & \mu_3^{(3)} + p_{31}^{(3)} m_{13} + p_{32}^{(3)} m_{23} \end{bmatrix}. \qquad (50)$$

Thus the process can be progressed from $M_{n-1}$ to $M_n$ using Theorems 4, 5 and 6.

## 5. Special case $N = 1$:

When the state reduction process results in a single state we in effect end up with a MRP $\{(X_k^{(1)}, T_k^{(1)}), k \geq 0\}$ on the state space $S_1 = \{1\}$. In this case the embedded irreducible MC $\{X_k^{(1)}, k \geq 0\}$ leads simply to $X_k^{(1)} \equiv 1$ for all $k$, having a single element transition matrix $P^{(1)} = [p_{11}^{(1)}] = [1]$. Thus the stationary probability distribution is $\pi_1^{(1)} = 1$.

Further the MRP reduces to the Renewal Process $\{T_k^{(1)}, k \geq 0\}$ where the distribution of the time between transitions, $Q_{11}^{(1)}(t) = F_{11}^{(t)}(t) = P\{T_{k+1}^{(1)} - T_k^{(1)} \leq t\}$. The mean state holding time $\mu_1^{(1)} = E[T_{k+1}^{(1)} - T_k^{(1)}]$. Since $\pi_1^{(1)} = 1$, the mean asymptotic increment $\lambda_1^{(1)} = \mu_1^{(1)}$ implying that the mean recurrence time is simply $m_{11} = \mu_1^{(1)}$.

## 6. Special case $N = 2$:

We consider the MRP $\{(X_k^{(2)}, T_k^{(2)}), k \geq 0\}$ on the state space $S_2 = \{1, 2\}$ with embedded irreducible MC $\{X_k^{(2)}, k \geq 0\}$ having a transition matrix $P^{(2)} = [p_{ij}^{(2)}]$ and mean state holding times $\mu_i^{(2)}$, $i = 1, 2$.

The state reduction procedure implies $\pi_1^{(2)} p_{21}^{(2)} = \pi_2^{(2)} p_{12}^{(2)}$, so that the stationary probabilities for the MC $\{X_k^{(2)}, k \geq 0\}$ are given by

$$\pi_1^{(2)} = \frac{p_{21}^{(2)}}{p_{12}^{(2)} + p_{21}^{(2)}}, \quad \pi_2^{(2)} = \frac{p_{12}^{(2)}}{p_{12}^{(2)} + p_{21}^{(2)}}. \tag{51}$$

For the $N = 2$ state situation, we have solved the matrix equation (22) when $n = 2$, in Section 4, in element form leading to equation (49) for $M_2$.

Note that from equation (37), with $n = 2$, $j = 1$, $S(2) = 1 - p_{22}^{(2)} = p_{21}^{(2)}$, so that $m_{21} = \frac{\mu_2^{(2)}}{S(2)}$, consistent with the expression for $m_{21}$ in equation (49).

Further, for $n = 2$, $i = 1$, equation (45) implies that $m_{12} = \frac{\nu_1^{(1,2)}}{1 - q_{11}^{(1,2)}} = \frac{\mu_1^{(2)}}{R(1,2)} = \frac{\mu_1^{(2)}}{1 - p_{11}^{(2)}} = \frac{\mu_1^{(2)}}{p_{12}^{(2)}}$, consistent with the result for $m_{12}$ in equation (49).

Note for the mean recurrence times, $m_{ii}$, we have from the proof of Theorem 3 that $m_{ii} = \frac{\lambda_1^{(2)}}{\pi_i^{(2)}}$ where $\lambda_1^{(2)} = \pi_1^{(2)} \mu_1^{(2)} + \pi_2^{(2)} \mu_2^{(2)}$.

The mean asymptotic increment is given by $\lambda_1^{(2)} = \pi_1^{(2)} \mu_1^{(2)} + \pi_2^{(2)} \mu_2^{(2)} = \frac{p_{21}^{(2)} \mu_1^{(2)} + p_{12}^{(2)} \mu_1^{(2)}}{p_{12}^{(2)} + p_{21}^{(2)}}$, implying $m_{11} = \frac{\lambda_1^{(2)}}{\pi_1^{(2)}} = \frac{p_{21}^{(2)} \mu_1^{(2)} + p_{12}^{(2)} \mu_2^{(2)}}{p_{12}^{(2)} + p_{21}^{(2)}} \cdot \frac{p_{12}^{(2)} + p_{21}^{(2)}}{p_{21}^{(2)}} = \frac{p_{21}^{(2)} \mu_1^{(2)} + p_{12}^{(2)} \mu_2^{(2)}}{p_{21}^{(2)}} = \mu_1^{(2)} + \frac{p_{12}^{(2)}}{p_{21}^{(2)}} \mu_2^{(2)} = \mu_1^{(1)}$ as already deduced for the $N = 1$ case.

Further, $m_{22} = \frac{\lambda_1^{(2)}}{\pi_2^{(2)}} = \frac{p_{21}^{(2)} \mu_1^{(2)} + p_{12}^{(2)} \mu_2^{(2)}}{p_{12}^{(2)}} = \frac{p_{21}^{(2)}}{p_{12}^{(2)}} \mu_1^{(2)} + \mu_2^{(2)}$, as given by equation (49).

When the MRP reduces to an irreducible MC, the stationary probabilities are as in equation (51), but the asymptotic mean increment is given by $\lambda_1^{(2)} = 1$, since $\mu_1^{(2)} = \mu_2^{(2)} = 1$, implying that

$$M_2 = \begin{bmatrix} m_{11} & m_{12} \\ m_{21} & m_{22} \end{bmatrix} = \begin{bmatrix} 1 + \frac{p_{12}^{(2)}}{p_{21}^{(2)}} & \frac{1}{p_{12}^{(2)}} \\ \frac{1}{p_{21}^{(2)}} & 1 + \frac{p_{21}^{(2)}}{p_{12}^{(2)}} \end{bmatrix},$$

as is well known (Hunter, Ex 7.3.4, [10]).

## 7. Special case $N = 3$

We examine the MRP $\{(X_k^{(3)}, T_k^{(3)}), k \geq 0\}$ on the state space $S_3 = \{1, 2, 3\}$ with embedded irreducible MC $\{X_k^{(3)}, k \geq 0\}$ having a transition matrix $P^{(3)} = [p_{ij}^{(3)}]$ and mean state holding times $\mu_i^{(3)}$, $i = 1, 2, 3$.

Firstly the stationary probabilities for the MC can be found from the state reduction process. From (15),

$$\pi^{(3)T} = (\pi_1^{(3)}, \pi_2^{(3)}, \pi_3^{(3)}) = c_2 (\pi_1^{(2)}, \pi_2^{(2)}, \frac{\sum_{i=1}^{2} \pi_i^{(2)} p_{i3}^{(3)}}{S(3)})$$

where $\pi_1^{(2)} = \dfrac{p_{21}^{(2)}}{p_{12}^{(2)} + p_{21}^{(2)}}$, $\pi_2^{(2)} = \dfrac{p_{12}^{(2)}}{p_{12}^{(2)} + p_{21}^{(2)}}$.

Let us introduce some notation that has previously been used in the literature that will simplify the expressions.

Define $\Delta_1 \equiv p_{21}^{(3)} p_{31}^{(3)} + p_{21}^{(3)} p_{32}^{(3)} + p_{23}^{(3)} p_{31}^{(3)}$, $\Delta_2 \equiv p_{12}^{(3)} p_{31}^{(3)} + p_{12}^{(3)} p_{32}^{(3)} + p_{13}^{(3)} p_{32}^{(3)}$,

$\Delta_3 \equiv p_{13}^{(3)} p_{21}^{(3)} + p_{12}^{(3)} p_{23}^{(3)} + p_{13}^{(3)} p_{23}^{(3)}$ and $\Delta \equiv \Delta_1 + \Delta_2 + \Delta_3$.

Now from (31) and (32) $p_{ij}^{(2)} = p_{ij}^{(3)} + \dfrac{p_{i3}^{(3)} p_{3j}^{(3)}}{S(3)}, (i,j) \in \{1,2\}$ where $S(3) = 1 - p_{33}^{(3)} = p_{31}^{(3)} + p_{32}^{(3)}$,

and $\mu_i^{(2)} = \mu_i^{(3)} + \dfrac{\mu_3^{(3)} p_{i3}^{(3)}}{S(3)}$, $(1 \leq i \leq 2)$, implying

$p_{12}^{(2)} = \dfrac{\Delta_2}{p_{31}^{(3)} + p_{32}^{(3)}}$, $p_{21}^{(2)} = \dfrac{\Delta_1}{p_{31}^{(3)} + p_{32}^{(3)}}$, $\mu_1^{(2)} = \mu_1^{(3)} + \dfrac{\mu_3^{(3)} p_{13}^{(3)}}{p_{31}^{(3)} + p_{32}^{(3)}}$ and $\mu_2^{(2)} = \mu_2^{(3)} + \dfrac{\mu_3^{(3)} p_{23}^{(3)}}{p_{31}^{(3)} + p_{32}^{(3)}}$.

Further

$$\pi_1^{(2)} = \dfrac{\Delta_1}{(p_{12}^{(3)} + p_{21}^{(3)})(p_{31}^{(3)} + p_{32}^{(3)}) + p_{13}^{(3)} p_{32}^{(3)} + p_{23}^{(3)} p_{31}^{(3)}} = c_2 \pi_1^{(3)},$$

$$\pi_2^{(2)} = \dfrac{\Delta_2}{(p_{12}^{(3)} + p_{21}^{(3)})(p_{31}^{(3)} + p_{32}^{(3)}) + p_{13}^{(3)} p_{32}^{(3)} + p_{23}^{(3)} p_{31}^{(3)}} = c_2 \pi_2^{(3)}.$$

$$\dfrac{\sum_{i=1}^{2} \pi_i^{(2)} p_{i3}^{(3)}}{S(3)} = \dfrac{\Delta_3}{(p_{12}^{(3)} + p_{21}^{(3)})(p_{31}^{(3)} + p_{32}^{(3)}) + p_{13}^{(3)} p_{32}^{(3)} + p_{23}^{(3)} p_{31}^{(3)}} = c_2 \pi_3^{(3)}.$$

implying that

$$\pi_i^{(3)} = \dfrac{\Delta_i}{\Delta}, i = 1, 2, 3. \tag{52}$$

Using the facts, derived from the above observations,

$R(2,3) = 1 - q_{22}^{(2,3)} = 1 - p_{22}^{(3)} = p_{21}^{(3)} + p_{23}^{(3)}$,

$R(1,3) = 1 - q_{11}^{(1,3)} = 1 - p_{11}^{(3)} - \dfrac{p_{12}^{(3)} p_{21}^{(3)}}{1 - p_{22}^{(3)}} = \dfrac{\Delta_3}{p_{21}^{(3)} + p_{23}^{(3)}}$.

$v_2^{(2,3)} = \mu_2^{(3)}$, $v_1^{(1,3)} = v_1^{(2,3)} + \dfrac{q_{12}^{(2,3)} v_2^{(2,3)}}{1 - q_{22}^{(2,3)}} = \mu_1^{(3)} + \dfrac{p_{12}^{(3)} \mu_2^{(3)}}{p_{21}^{(3)} + p_{23}^{(3)}}$,

and using the simplifications that $\Delta_1 p_{13}^{(3)} + \Delta_2 p_{23}^{(3)} = \Delta_3 (p_{31}^{(3)} + p_{32}^{(3)})$,

$p_{23}^{(3)} p_{32}^{(3)} + \Delta_1 = (p_{21}^{(3)} + p_{23}^{(3)})(p_{31}^{(3)} + p_{32}^{(3)})$, $p_{12}^{(3)} p_{31}^{(3)} + \Delta_2 = (p_{12}^{(3)} + p_{13}^{(3)})(p_{31}^{(3)} + p_{32}^{(3)})$ and

$p_{12}^{(3)} p_{21}^{(3)} + \Delta_3 = (p_{12}^{(3)} + p_{13}^{(3)})(p_{21}^{(3)} + p_{23}^{(3)})$ we express all the elemental expressions of the $M_3$ matrix for the MFPTs in terms of the $p_{ij}^{(3)}$ and the $\mu_i^{(3)}$. This leads to

$$M_3 = \begin{bmatrix} \dfrac{\Delta_1\mu_1^{(3)} + \Delta_2\mu_2^{(3)} + \Delta_3\mu_3^{(3)}}{\Delta_1} & \dfrac{(p_{31}^{(3)} + p_{32}^{(3)})\mu_1^{(3)} + p_{13}^{(3)}\mu_3^{(3)}}{\Delta_2} & \dfrac{(p_{21}^{(3)} + p_{23}^{(3)})\mu_1^{(3)} + p_{12}^{(3)}\mu_2^{(3)}}{\Delta_3} \\ \dfrac{(p_{31}^{(3)} + p_{32}^{(3)})\mu_2^{(3)} + p_{23}^{(3)}\mu_3^{(3)}}{\Delta_1}, & \dfrac{\Delta_1\mu_1^{(3)} + \Delta_2\mu_2^{(3)} + \Delta_3\mu_3^{(3)}}{\Delta_2} & \dfrac{p_{21}^{(3)}\mu_1^{(3)} + (p_{12}^{(3)} + p_{13}^{(3)})\mu_2^{(3)}}{\Delta_3} \\ \dfrac{p_{32}^{(3)}\mu_2^{(3)} + (p_{21}^{(3)} + p_{23}^{(3)})\mu_3^{(3)}}{\Delta_1} & \dfrac{p_{31}^{(3)}\mu_1^{(3)} + (p_{12}^{(3)} + p_{13}^{(3)})\mu_3^{(3)}}{\Delta_2} & \dfrac{\Delta_1\mu_1^{(3)} + \Delta_2\mu_2^{(3)} + \Delta_3\mu_3^{(3)}}{\Delta_3} \end{bmatrix}. \quad (53)$$

Note that the expression for $m_{33}$ can also be deduced from either equation (34 or (36). Note also that from the properties of MRPs, $m_{ii} = \dfrac{\lambda_1^{(3)}}{\pi_i^{(3)}}$. The diagonal elements of (53) are consistent with this observation since, using (52), the mean asymptotic increment is given by

$$\lambda_1^{(3)} = \pi_1^{(3)}\mu_1^{(3)} + \pi_2^{(3)}\mu_2^{(3)} + \pi_3^{(3)}\mu_3^{(3)} = \dfrac{\Delta_1\mu_1^{(3)} + \Delta_2\mu_2^{(3)} + \Delta_3\mu_3^{(3)}}{\Delta}.$$

For the MC case, $\mu_i^{(3)} = 1$, $i = 1, 2, 3$. Substituting and simplifying (53) yields

$$M_3 = \begin{bmatrix} \dfrac{\Delta}{\Delta_1} & \dfrac{p_{13}^{(3)} + p_{31}^{(3)} + p_{32}^{(3)}}{\Delta_2} & \dfrac{p_{12}^{(3)} + p_{21}^{(3)} + p_{23}^{(3)}}{\Delta_3} \\ \dfrac{p_{23}^{(3)} + p_{31}^{(3)} + p_{32}^{(3)}}{\Delta_1} & \dfrac{\Delta}{\Delta_2} & \dfrac{p_{12}^{(3)} + p_{13}^{(3)} + p_{21}^{(3)}}{\Delta_3} \\ \dfrac{p_{21}^{(3)} + p_{23}^{(3)} + p_{32}^{(3)}}{\Delta_1} & \dfrac{p_{12}^{(3)} + p_{13}^{(3)} + p_{31}^{(3)}}{\Delta_2} & \dfrac{\Delta}{\Delta_3} \end{bmatrix}. \quad (54)$$

These results are equivalent to those given in Example 3.2 of [11] where it is shown that the MC, with the transition matrix $P^{(3)} = \left[p_{ij}^{(3)}\right]$, is irreducible (and hence a stationary distribution exists) if and only if $\Delta_1 > 0$, $\Delta_2 > 0$, $\Delta_3 > 0$ with stationary distribution given as in (52) and MFPT matrix is given by (54).

## 8. Special case $N = 4$

We examine the MRP $\{(X_k^{(4)}, T_k^{(4)}), k \geq 0\}$ on the state space $S_4 = \{1,2,3,4\}$ with embedded irreducible MC $\{X_k^{(4)}, k \geq 0\}$ having a transition matrix $P^{(4)} = \left[p_{ij}^{(4)}\right]$ and mean state holding times $\mu_i^{(4)}$, $i = 1, 2, 3, 4$.

We extend the $M_3$ matrix, using Theorem 5, equation (37) for the $m_{4j}$ terms, Theorem 6, equation (47) for the $m_{i4}$ terms and Theorem 4, equation (36) for the $m_{44}$. This leads to the pattern of the calculations that need to be carried out in the boundary column and row.

$$M_4 = \begin{bmatrix} m_{11} & m_{12} & m_{13} & \dfrac{v_1^{(1,4)}}{R(1,4)} \\ m_{21} & m_{22} & m_{23} & \dfrac{v_2^{(2,4)} + q_{21}^{(2,4)} m_{14}}{R(2,4)} \\ m_{31} & m_{32} & m_{33} & \dfrac{v_3^{(3,4)} + q_{31}^{(3,4)} m_{14} + q_{32}^{(3,4)} m_{24}}{R(3,4)} \\ \dfrac{\mu_4^{(4)} + p_{42}^{(4)} m_{21} + p_{43}^{(4)} m_{31}}{S(4)} & \dfrac{\mu_4^{(4)} + p_{41}^{(4)} m_{12} + p_{43}^{(4)} m_{32}}{S(4)} & \dfrac{\mu_4^{(4)} + p_{41}^{(4)} m_{13} + p_{42}^{(4)} m_{23}}{S(4)} & \mu_4^{(4)} + p_{41}^{(4)} m_{14} + p_{42}^{(4)} m_{24} + p_{43}^{(4)} m_{34} \end{bmatrix}.$$

The determination of the terms in the fourth column require careful computation.

Firstly we express the required MFPTs in terms of the initial transition probabilities, the $p_{ij}^{(4)}$ computing successively:

$$q_{ij}^{(3,4)} = p_{ij}^{(4)}, \quad i=1,2,3, \; j=1,2,3,$$

$$q_{ij}^{(2,4)} = q_{ij}^{(3,4)} + \frac{q_{i3}^{(3,4)} q_{3j}^{(3,4)}}{1 - q_{33}^{(3,4)}} = p_{ij}^{(4)} + \frac{p_{i3}^{(4)} p_{3j}^{(4)}}{1 - p_{33}^{(4)}}, \quad i=1,2, \; j=1,2,$$

$$q_{11}^{(1,4)} = q_{11}^{(2,4)} + \frac{q_{12}^{(2,4)} q_{21}^{(2,4)}}{1 - q_{22}^{(2,4)}} = \left( p_{11}^{(4)} + \frac{p_{13}^{(4)} p_{31}^{(4)}}{1 - p_{33}^{(4)}} \right) + \frac{\left( p_{12}^{(4)} + \dfrac{p_{13}^{(4)} p_{32}^{(4)}}{1 - p_{33}^{(4)}} \right)\left( p_{21}^{(4)} + \dfrac{p_{23}^{(4)} p_{31}^{(4)}}{1 - p_{33}^{(4)}} \right)}{1 - \left( p_{22}^{(4)} + \dfrac{p_{23}^{(4)} p_{32}^{(4)}}{1 - p_{33}^{(4)}} \right)}.$$

Further $R(i,4) = 1 - q_{ii}^{(i,4)}$, $i=1,2,3$ so that

$$R(3,4) = 1 - p_{33}^{(4)} = p_{31}^{(4)} + p_{32}^{(4)} + p_{34}^{(4)},$$

$$R(2,4) = 1 - q_{22}^{(2,4)} = 1 - p_{22}^{(4)} - \frac{p_{23}^{(4)} p_{32}^{(4)}}{1 - p_{33}^{(4)}}$$

$$= \frac{(p_{21}^{(4)} + p_{24}^{(4)})(p_{31}^{(4)} + p_{34}^{(4)}) + p_{23}^{(4)}(p_{31}^{(4)} + p_{34}^{(4)}) + (p_{21}^{(4)} + p_{24}^{(4)}) p_{32}^{(4)}}{p_{31}^{(4)} + p_{32}^{(4)} + p_{34}^{(4)}}.$$

$$R(1,4) = 1 - q_{11}^{(1,4)} = 1 - q_{11}^{(2,4)} - \frac{q_{12}^{(2,4)} q_{21}^{(2,4)}}{1 - q_{22}^{(2,4)}}$$

$$= 1 - \left( p_{11}^{(4)} - \frac{p_{13}^{(4)} p_{31}^{(4)}}{1 - p_{33}^{(4)}} \right) - \frac{\left( p_{12}^{(4)} - \dfrac{p_{13}^{(4)} p_{32}^{(4)}}{1 - p_{33}^{(4)}} \right)\left( p_{21}^{(4)} - \dfrac{p_{23}^{(4)} p_{31}^{(4)}}{1 - p_{33}^{(4)}} \right)}{1 - \left( p_{22}^{(4)} - \dfrac{p_{23}^{(4)} p_{32}^{(4)}}{1 - p_{33}^{(4)}} \right)}.$$

Also $v_i^{(3,4)} = \mu_i^{(4)}$, $i = 1, 2, 3$, and, from equation (43),

$$v_i^{(2,4)} = v_i^{(3,4)} + \frac{q_{i3}^{(3,4)} v_3^{(3,4)}}{1 - q_{33}^{(3,4)}} = \mu_i^{(4)} + \frac{p_{i3}^{(4)} \mu_3^{(4)}}{1 - p_{33}^{(4)}}, \quad i=1,2.$$

$$v_1^{(1,4)} = v_1^{(2,4)} + \frac{q_{12}^{(2,4)} v_2^{(2,4)}}{1 - q_{22}^{(2,4)}} = \left( \mu_1^{(4)} + \frac{p_{13}^{(4)} \mu_3^{(4)}}{1 - p_{33}^{(4)}} \right) + \frac{\left( p_{12}^{(4)} + \dfrac{p_{13}^{(4)} p_{32}^{(4)}}{1 - p_{33}^{(4)}} \right)\left( \mu_2^{(4)} + \dfrac{p_{23}^{(4)} \mu_3^{(4)}}{1 - p_{33}^{(4)}} \right)}{1 - \left( p_{22}^{(4)} + \dfrac{p_{23}^{(4)} p_{32}^{(4)}}{1 - p_{33}^{(4)}} \right)}.$$

We obtain, $m_{14} = \dfrac{v_1^{(2,4)} + \left[\dfrac{q_{12}^{(2,4)}v_2^{(2,4)}}{1-q_{22}^{(2,4)}}\right]}{1-q_{11}^{(2,4)} - \dfrac{q_{12}^{(2,4)}q_{21}^{(2,4)}}{1-q_{22}^{(2,4)}}} = \dfrac{(1-q_{22}^{(2,4)})v_1^{(2,4)} + q_{12}^{(2,4)}v_2^{(2,4)}}{(1-q_{11}^{(2,4)})(1-q_{22}^{(2,4)}) - q_{12}^{(2,4)}q_{21}^{(2,4)}}.$

$m_{24} = \dfrac{\left[p_{21}^{(4)}(1-p_{33}^{(4)}) + p_{23}^{(4)}p_{31}^{(4)}\right]m_{14} + \left[\mu_2^{(4)}(1-p_{33}^{(4)}) + \mu_3^{(4)}p_{23}^{(4)}\right]}{(1-p_{22}^{(4)})(1-p_{33}^{(4)}) - p_{23}^{(4)}p_{32}^{(4)}}.$

$m_{34} = \dfrac{p_{31}^{(4)}m_{14} + p_{32}^{(4)}m_{24} + \mu_3^{(4)}}{p_{31}^{(4)} + p_{32}^{(4)} + p_{34}^{(4)}}.$

These expressions can be simplified upon substitution of the terms above, but the numerators and denominators are not particularly simple expressions.

$m_{14} = \dfrac{N_{14}}{D_{14}}$ where

$N_{14} = \left(1 - p_{22}^{(4)} - \dfrac{p_{23}^{(4)}p_{32}^{(4)}}{1-p_{33}^{(4)}}\right)\left(\mu_1^{(4)} + \dfrac{p_{13}^{(4)}\mu_3^{(4)}}{1-p_{33}^{(4)}}\right) + \left(p_{12}^{(4)} + \dfrac{p_{13}^{(4)}p_{32}^{(4)}}{1-p_{33}^{(4)}}\right)\left(\mu_2^{(4)} + \dfrac{p_{23}^{(4)}\mu_3^{(4)}}{1-p_{33}^{(4)}}\right)$

$D_{14} = \left(1 - p_{11}^{(4)} - \dfrac{p_{13}^{(4)}p_{31}^{(4)}}{1-p_{33}^{(4)}}\right)\left(1 - p_{22}^{(4)} - \dfrac{p_{23}^{(4)}p_{32}^{(4)}}{1-p_{33}^{(4)}}\right) + \left(p_{12}^{(4)} + \dfrac{p_{13}^{(4)}p_{32}^{(4)}}{1-p_{33}^{(4)}}\right)\left(p_{21}^{(4)} + \dfrac{q_{23}^{(4)}q_{31}^{(4)}}{1-p_{33}^{(4)}}\right).$

## 9. Algorithms for the computation of the matrix of MFPTs

From the special cases considered in Sections 5 to 8, we have typically extended the calculations for the elements of $M_n$ from the matrix $M_{n-1}$ by appending the elements $m_{in}$ ($i = 1, .., n-1$), $m_{nn}$ and $m_{nj}$ ($j = 1, \ldots, n-1$). However by exploring these calculations in more depth it is apparent that an recursive process for computing the $m_{ij}$ elements for the original $M_N$ matrix can be constructed through three separate different procedures corresponding to the three cases $i < j$, $i = j$, and $i > j$. We can separate these three cases as follows, using the notation developed earlier, viz., $p_{ij}^{(n)}, \mu_i^{(n)}, S(n)$, as given by equations (31) and (32), $q_{ij}^{(t-1,n)}$ with $q_{ij}^{(n-1,n)} = p_{ij}^{(n)}$ as given by equations (41) and (42), $v_i^{(t-1,n)}$ with $v_i^{(n-1,n)} = \mu_i^{(n)}$ as given by equations (43) and (44), and $R(i,n) = 1 - q_{ii}^{(i,n)}$, as given by equation (46).

**Theorem 7**: The elements of the MFPT matrix $M_N = [m_{ij}]$ ($i = 1,\ldots,N, j = 1,2,\ldots,N$) for the MRP $\{(X_k^{(N)}, T_k^{(N)}), k \geq 0\}$ with state space $S_N = \{1, 2, \ldots, N\}$, transition matrix $P^{(N)} = [p_{ij}^{(N)}]$ and vector of mean holding times $\boldsymbol{\mu}^{(N)T} = (\mu_1^{(N)}, \ldots, \mu_N^{(N)})$ can be expressed as follows, where we use the notation

(a) $m_{ij} = \dfrac{\mu_i^{(i)} + \sum_{k=1, k\neq j}^{i-1} p_{ik}^{(i)}m_{kj}}{S(i)}$, ($i = 3,\ldots,N; j = 1,\ldots,i-1$), (55)

with $m_{21} = \dfrac{\mu_2^{(2)}}{S(2)}$. (56)

(b) $m_{ii} = \mu_i^{(i)} + \sum_{k=1}^{i-1} p_{ik}^{(i)}m_{ki}$, ($i = 2,\ldots,N$), (57)

with $m_{11} = \mu_1^{(1)}$. (58)

(c) $$m_{ij} = \frac{v_i^{(i,j)} + \sum_{k=1}^{i-1} q_{ik}^{(i,j)} m_{kj}}{R(i,j)}, \quad (i = 2, 3, ..., N-1; \; j = i+1, ....., N),$$ (59)

with $m_{1j} = \dfrac{v_1^{(1,j)}}{R(1,j)}, \; (j = 2, 3, ..., N).$ (60)

**Proof:** The expressions in (a) are when the indices $i > j$, in (b) when $i = j$ and in (c) when $i < j$. Equation (55) follows from (37) and (56) from (49); equation (57) from (36) and (58) from (36); equation (59) and (60) from (47).

□

The general procedure described by Theorem 7 is difficult to program, for a general state space, using MatLab. In particular the computation for the $m_{ij}$ when $j > i$ demands additional computations. Typically the elements of $P^{(n)} = \left[ p_{ij}^{(n)} \right]$, an $n \times n$ stochastic matrix, are easily found by the GTH algorithm. However in order to compute the $q_{ij}^{(t,n)}$ requires first identifying the starting elements of $Q_n^{(n-1)} = \left[ q_{ij}^{(n-1,n)} \right]_{(n-1)\times(n-1)} = \left[ p_{ij}^{(n)} \right]_{(n-1)\times(n-1)}$ i.e. the elements are from the sub-stochastic matrix found from the first $n - 1$ rows and $n - 1$ columns of $P^{(n)}$. The reduction through the sequence of GTH reduction procedures leads from $Q_n^{(n-1)} = \left[ q_{ij}^{(n-1,n)} \right]_{(n-1)\times(n-1)} \to Q_n^{(n-2)} = \left[ q_{ij}^{(n-2,n)} \right]_{(n-2)\times(n-2)} \to ....$ to eventually arrive at $Q_n^{(2)} = \left[ q_{ij}^{(2,n)} \right]_{2\times 2}$ and finally at $Q_n^{(1)} = \left[ q_{ij}^{(1,n)} \right]_{1\times 1}$. Because of the truncation of $P^{(n)}$, for each $n$, to start with $Q_n^{(n-1)}$, this process has to be implemented for each value of $n = N, N-1, ...., 2$. Thus the GTH algorithmic reduction has to be carried out a number of times, as in the following grid,

$P^{(N)} \to Q_{N-1}^{(N)} \to .... \to Q_{21}^{(N)} \to Q_1^{(N)}$ followed by $P^{(N-1)} \to Q_{N-2}^{(N-1)} \to Q_{N-2}^{(N-1)} \to .... \to Q_1^{(N-1)}$, leading successively to $P^{(3)} \to Q_2^{(3)} \to Q_1^{(3)}$ and finally to $P^{(2)} \to Q_1^{(2)}$. From the initial $(N-1)$ GTH algorithmic procedures for the $p_{ij}^{(n)}$, followed by $(N-1)$ matrix reductions to start with the initial $q_{ij}^{(n-1,n)}$ there are a further $(N-2) + (N-3) + ....+ 1$ reductions for a total $N(N-1)/2$ separate GTH procedures. For the original GTH procedure for finding the stationary probabilities we only needed retention of the $p_{in}^{(n)}$ and $p_{nj}^{(n)}$ boundary terms of the $P^{(n)}$ matrices, whereas for the MFPT's we need to retain additional elements of the $P^{(n)}$ leading to the $Q_i^{(n)}$ matrices.

In the computation of the $R(i,n) = 1 - q_{ii}^{(1,n)}$ expressions we have no easy technique to ensure that no subtraction is required. This is due to the fact that the sum of the elements in the last row of the $Q_i^{(n)}$ matrix do not sum to 1, as in the $P^{(n)}$ matrices.

The $q_{ik}^{(t,n)}$ terms only arise in the computation of the MFPTs $m_{ij}$ when $i < j$, whereas the $p_{ij}^{(n)}$ are all that is needed to compute the $m_{ij}$ when $i > j$ and these probabilities are all that is needed to compute the mean holding times $\mu_i^{(n)}$.

The above observations lead to the following as a general technique for finding all the elements of $M$ for the case of a given MRP. Since we effectively use the computations of the GTH procedure, we call this the "Extended GTH" (EGTH) algorithm.

**EGTH Algorithm**

Step 1(*i*): Start with $P^{(N)} = \left[p_{ij}^{(N)}\right]$, carry out the GTH algorithm by calculating successively,

for $n = N, N-1, \ldots, 2, p_{ij}^{(n-1)} = p_{ij}^{(n)} + \dfrac{p_{in}^{(n)} p_{nj}^{(n)}}{S(n)}$, $1 \le i \le n-1, 1 \le j \le n-1$, where $S(n) = \sum_{j=1}^{n-1} p_{nj}^{(n)}$.

(Note that we only have to retain the $p_{in}^{(n)} (1 \le i \le n-1)$ and $p_{nj}^{(n)} (1 \le j \le n-1)$, i.e. the *n*-th row and *n*-th column of $P^{(n)}$ for $n = 2, \ldots, N$, as in the GTH algorithm.)

Step 1(*ii*): Start with the mean holding time vector $\boldsymbol{\mu}^{(N)T} = (\mu_1^{(N)}, \mu_2^{(N)}, \ldots, \mu_{N-1}^{(N)}, \mu_N^{(N)})$ and calculate successively for $n = N, N-1, \ldots, 2$, $\mu_i^{(n-1)} = \mu_i^{(n)} + \dfrac{\mu_n^{(n)} p_{in}^{(n)}}{S(n)}$, $1 \le i \le n-1$.

Step 1(*iii*): Calculate the $N \times 1$ column vector $\boldsymbol{m}_N^{(1)(N)} = (m_{i1})$, where $m_{11} = \mu_1^{(1)}$,

$m_{21} = \dfrac{\mu_2^{(2)}}{S(2)}$, and for $i = 3, \ldots, N$, $m_{i1} = \dfrac{\mu_i^{(i)} + \sum_{k=2}^{i-1} p_{ik}^{(i)} m_{k1}}{S(i)}$.

This gives the entries of the first column of $M = [m_{ij}]$, i.e. $\boldsymbol{m}_N^{(1)(N)}$ where $M = \left(\boldsymbol{m}_N^{(1)(N)}, \boldsymbol{m}_N^{(2)(N)}, \ldots, \boldsymbol{m}_N^{(N)(N)}\right)$ with $\boldsymbol{m}_N^{(1)(N)T} = (m_{11}, m_{21}, \ldots, m_{N1})$.

Comment: The steps that follow are based on the observation that by starting with $P^{(N)}$, which we define as $P^{(N)(1)}$, we are able to find expressions for $\boldsymbol{m}_N^{(1)(N)}$, the first column of *M*, giving the MFPTs to state 1 from all the other states. Successively we permute the elements of $P^{(N)}$ so as to do this for each of the states 2, …,*N*. For state 2 we can do this by moving the elements of first column of $P^{(N)}$ to after the *N*-th column, followed by moving the first row to the last row, to obtain a new transition matrix $P^{(N)(2)}$.

$$P^{(N)} \equiv P^{(N)(1)} = \begin{bmatrix} p_{11} & p_{12} & \cdots & p_{1,N-1} & p_{1,N} \\ p_{21} & p_{22} & \cdots & & \\ \cdot\cdot & \cdot\cdot & \cdots & \cdot\cdot & \cdot\cdot \\ p_{N-1,1} & p_{N-1,2} & \cdots & p_{N-1,N-1} & p_{N-1,N} \\ p_{N1} & p_{N2} & \cdots & p_{N,N-1} & p_{NN} \end{bmatrix}$$

$$\rightarrow \begin{bmatrix} p_{12} & \cdots & p_{1,N-1} & p_{1,N} & p_{11} \\ p_{22} & \cdots & p_{2,N-1} & p_{2N} & p_{21} \\ \cdot\cdot & \cdots & \cdots & \cdot\cdot & \cdot\cdot \\ p_{N-1,2} & \cdots & p_{N-1,N-1} & p_{N-1,N} & p_{N-1,1} \\ p_{N2} & \cdots & p_{N,N-1} & p_{N,N} & p_{N1} \end{bmatrix} \rightarrow \begin{bmatrix} p_{22} & \cdots & p_{2,N-1} & p_{2N} & p_{21} \\ \cdot\cdot & \cdots & \cdot\cdot & \cdot\cdot & \cdot\cdot \\ p_{N-1,2} & \cdots & p_{N-1,N-1} & p_{N-1,N} & p_{N-1,1} \\ p_{N2} & \cdots & p_{N,N-1} & p_{NN} & p_{N1} \\ p_{12} & \cdots & p_{1,N-1} & p_{1N} & p_{11} \end{bmatrix} \equiv P^{(N)(2)}$$

There are a variety of ways we can do this. Here are three such ways:

(*i*) Let $\boldsymbol{e}_i^T = (0,0,\ldots,1,0,\ldots,0)$ be the *i*-th elementary row vector with 1 in the *i*-th position and 0 elsewhere and $\boldsymbol{e}_i$ is the *i*-th elementary column vector.
Let $R_1^{(N)} = [\boldsymbol{e}_N, \boldsymbol{e}_1, \ldots, \boldsymbol{e}_{N-1}]$ and $C_1^{(N)} = [\boldsymbol{e}_2, \ldots, \boldsymbol{e}_N, \boldsymbol{e}_1]$. Then $P^{(N)(2)} = R_1^{(N)} P^{(N)(1)} C_1^{(N)}$.

(*ii*) $P^{(N)(2)} (mod(row + N - 2, N) + 1, mod(col + N - 2, N) + 1) = P^{(N)(1)} (row, col)$.

This can be done in stages if necessary with say the row shift followed by the column shift.

(*iii*) In MatLab use the "circshift" operator with $P^{(N)(2)} = circshift(P^{(N)(1)}, [-1,-1])$.

Step 2: For $k = 2, 3, 4, \ldots, N-1, N$.
(*i*) Repeat Step 1(*i*) but with $P^{(N)} = P^{(N)(k)}$ where $P^{(N)(k)} = R_1^{(N)} P^{(N)(k-1)} C_1^{(N)}$ with $P^{(N)(1)} = P^{(N)}$
(Comment: This steps leads to the appropriate $p_{in}^{(n)}$ and $p_{nj}^{(n)}$ elements.)
(*ii*) Repeat Step 1(*ii*) but with $\boldsymbol{\mu}^{(N)} = \boldsymbol{\mu}^{(N)(k)}$ where $\boldsymbol{\mu}^{(N)(k)T} = \boldsymbol{\mu}^{(N)(k-1)T} C_1^{(N)}$ with $\boldsymbol{\mu}^{(N)(1)} = \boldsymbol{\mu}^{(N)}$
(Comment: This step leads to the appropriate $\mu_i^{(n)}$ elements. In the case of a MC no permutation of the elements is required, since $\mu_i^{(N)} = 1$ for all *i*.)
(*iii*) Repeat Step 1(*iii*) to calculate the $N \times 1$ column vector $\bar{\boldsymbol{m}}_N^{(k)(N)}$ where
$\bar{\boldsymbol{m}}_N^{(k)(N)T} = (m_{kk}, m_{k+1,k}, \ldots, m_{Nk}, m_{1k}, \ldots, m_{k-1,k})$ .

Step 3: Combine the results of the Steps 1(*iii*) and 2(*iii*) to find *M* as follows.
Let $\bar{M} = (\bar{\boldsymbol{m}}_N^{(1)(N)}, \bar{\boldsymbol{m}}_N^{(2)(N)}, \ldots, \bar{\boldsymbol{m}}_N^{(N)(N)})$ and reorder the elements of $\bar{M}$ to obtain $M = (\boldsymbol{m}_N^{(1)(N)}, \boldsymbol{m}_N^{(2)(N)}, \ldots, \boldsymbol{m}_N^{(N)(N)})$. This can be carried out in MatLab by noting that for each row and column entry, $\bar{M}(\mod(row + col - 2, N) + 1, col) = M(row, col)$.

While this EGTH procedure requires *N* separate repetitions, one would have to carry *N* auxiliary sets of calculations to determine the $v_i^{(t,n)}$, as in Theorem 7, as well as retaining more calculations enroute to the derivation of the matrix *M*.

Another key observation is that the EGTH algorithm, as outlined above, retains the calculation accuracy with no subtractions being involved.

Note also that we do not compute the stationary probabilities in determining the MFPTs. In the MC setting they would typically be found using the basic GTH algorithm. However in this setting the stationary probabilities can also be found directly as the inverse of the $m_{ii}$, alleviating the necessity of any prior calculation. For example in the MC setting, the initial holding times are $\mu_i^{(N)} = 1$, we have from the first step in the EGTH algorithm that $m_{11} = \mu_1^{(1)}$ giving an alternative derivation of $\pi_1$ as $\pi_1 = 1/\mu_1^{(1)}$. Once again, no subtraction operation need be performed.

**10.  The Test Problems**

The following test problems were introduced by Harrod & Plemmons ([5]) and have been considered by others in different contexts. They were initially introduced as poorly conditioned examples for computing the stationary distribution of the underlying irreducible MC. While the dimensions of the state space are relatively small, the test problems lead to some computational difficulties.

**TP1:** (As modified by Heyman and Reeves ([8])

$$\begin{bmatrix} .1 & .6 & 0 & .3 & 0 & 0 \\ .5 & .5 & 0 & 0 & 0 & 0 \\ .5 & .2 & 0 & 0 & .3 & 0 \\ 0 & .7 & 0 & .2 & 0 & .1 \\ .1 & 0 & .8 & 0 & 0 & .1 \\ .4 & 0 & .4 & 0 & 0 & .2 \end{bmatrix}$$

**TP2**: (See also Benzi ([1]))

$$\begin{bmatrix} .85 & 0 & .149 & .0009 & 0 & .00005 & 0 & .00005 \\ .1 & .65 & .249 & 0 & .0009 & .00005 & 0 & .00005 \\ .1 & .8 & .09996 & .0003 & 0 & 0 & .0001 & 0 \\ 0 & .0004 & 0 & .7 & .2995 & 0 & .0001 & 0 \\ .0005 & 0 & .0004 & .399 & .6 & .0001 & 0 & 0 \\ 0 & .00005 & 0 & 0 & .00005 & .6 & .2499 & .15 \\ .00003 & 0 & .00003 & .00004 & 0 & .1 & .8 & .0999 \\ 0 & .00005 & 0 & 0 & .00005 & .1999 & .25 & .55 \end{bmatrix}.$$

**TP3**:

$$\begin{bmatrix} 0.999999 & 1.0\,E-07 & 2.0\,E-07 & 3.0\,E-07 & 4.0\,E-07 \\ 0.4 & 0.3 & 0 & 0 & 0.3 \\ 5.0\,E-07 & 0 & 0.999999 & 0 & 5.0\,E-07 \\ 5.0\,E-07 & 0 & 0 & 0.999999 & 5.0\,E-07 \\ 2.0\,E-07 & 3.0\,E-07 & 1.0\,E-07 & 4.0\,E-07 & 0.999999 \end{bmatrix}.$$

**TP4 and variants:**
**TP41**: $\varepsilon = 1.0E-01$; **TP42**: $\varepsilon = 1.0E-03$; **TP43**: $\varepsilon = 1.0E-05$; **TP44**: $\varepsilon = 1.0E-07$

$$\begin{bmatrix} .1-\varepsilon & .3 & .1 & .2 & .3 & \varepsilon & 0 & 0 & 0 & 0 \\ .2 & .1 & .1 & .2 & .4 & 0 & 0 & 0 & 0 & 0 \\ .1 & .2 & .2 & .4 & .1 & 0 & 0 & 0 & 0 & 0 \\ .4 & .2 & .1 & .2 & .1 & 0 & 0 & 0 & 0 & 0 \\ .6 & .3 & 0 & 0 & .1 & 0 & 0 & 0 & 0 & 0 \\ \varepsilon & 0 & 0 & 0 & 0 & .1-\varepsilon & .2 & .2 & .4 & .1 \\ 0 & 0 & 0 & 0 & 0 & .2 & .2 & .1 & .3 & .2 \\ 0 & 0 & 0 & 0 & 0 & .1 & .5 & 0 & .2 & .2 \\ 0 & 0 & 0 & 0 & 0 & .5 & .2 & .1 & 0 & .2 \\ 0 & 0 & 0 & 0 & 0 & .1 & .2 & .2 & .3 & .2 \end{bmatrix}$$

We carry out all the calculations using the academic version of MatLab (R2015b, 64bit on a MacBook Air). We first calculate the MFPT matrix $M$ for each of the given test problems, using the EGTH algorithm, under double precision. See Appendix 1 for the relevant MatLab code. (In Appendix 2, which appears only in the arXiv.com version of this paper, we present the accurate calculations for all the MFPTs for these test problems, as such results do not appear in the literature. We also give expressions for the relevant stationary probability vectors.)

We compute, for each test problem, with specified transition matrix, the following errors for the MFPT matrix, $M = [m_{ij}]$, given by the EGTH calculation, under both double and single precision: *Minimum absolute error* $= \min_{1 \leq i \leq m, 1 \leq j \leq m} \left| m_{ij} - \sum_{k \neq j} p_{ik} m_{kj} - 1 \right|$, *Maximum absolute error* $= \max_{1 \leq i \leq m, 1 \leq j \leq m} \left| m_{ij} - \sum_{k \neq j} p_{ik} m_{kj} - 1 \right|$, and the *overall residual error* $= \sum_{i=1}^{m} \sum_{j=1}^{m} \left| m_{ij} - \sum_{k \neq j} p_{ik} m_{kj} - 1 \right|$. These errors are given in Tables 1 and 2 below.

Table 1: Errors for the MFPTs (Double Precision)

| Test Problem | Min Abs Error | Max Absolute Error | Residual Error |
|---|---|---|---|
| TP1 | 0 | 1.1369E-13 | 2.9177E-13 |
| TP2 | 0 | 3.6380E-12 | 2.7776E-11 |
| TP3 | 0 | 1.8626E-09 | 2.7940E-09 |
| TP41 | 0 | 1.4211E-14 | 2.7337E-13 |
| TP42 | 0 | 1.8190E-12 | 1.9142E-11 |
| TP43 | 0 | 1.1642E-10 | 1.5717E-09 |
| TP44 | 0 | 7.4506E-09 | 1.4156E-07 |

Table 2: Errors for the MFPTs (Single Precision)

| Test Problem | Min Abs Error | Max Absolute Error | Residual Error |
|---|---|---|---|
| TP1 | 0 | 6.1035E-05 | 1.6773E-04 |
| TP2 | 0 | 1.9531E-03 | 1.3889E-02 |
| TP3 | 0 | 0.5000 | 3.0757 |
| TP41 | 0 | 7.6294E-06 | 1.0628E-04 |
| TP42 | 0 | 4.8828E-04 | 0.0050 |
| TP43 | 0 | 0.0860 | 0.7809 |
| TP44 | 0 | 5 | 85.8835 |

The test problems have been used as examples for testing various different algorithms for computing $M$, the matrix of MFPTs. In particular Heyman and O'Leary ([7]) considered five different procedures for computing the fundamental matrix $Z$, the group inverses $A^{\#}$ and $M$ (since these are all interconnected). Further Heyman and Reeves ([8]) also considered four different techniques for $M$ with their most accurate procedure based upon a state reduction procedure. We do not go into details of the procedures that they used but they compared the accuracy of the procedures by evaluating the number of accurate digits. The most accurate procedure in [7] was based upon using an LU factorization and normalization related to a state reduction procedure. In [8] the comparable procedure was also a state reduction procedure. The double precision result was used as the "true" result and the single precision result as the "computed" result. The number of accurate digits was defined as the overall average of $-\log_{10} \left| \frac{result_{true} - result_{computed}}{result_{true}} \right|$. Each of these two papers presented the results in figures and no actual numerical results were tabulated. We computed this statistic for each the seven test problems achieving the following results:

Table3: Average number of accurate digits

| TP 1  | 7.3504* |
|-------|---------|
| TP 2  | 7.2928  |
| TP 3  | 7.3526  |
| TP 41 | 7.3681  |
| TP 42 | 7.4157  |
| TP 43 | 7.4296  |
| TP 44 | 7.3321  |

*Note that for TP1 the MFPT from state 2 to state 1, is 2.00 for both the accurate and computed results so that the accurate digit quantity is infinite. The average in this case is taken over the remaining 35 possible pairs of states.

Considering the results of Heyman and O'Leary [7] and Heyman and Reeves [8], it is obvious that no procedure that they considered has any improvement over the results of this paper. Heyman and O'Leary have values between 6 and 7 for their favoured algorithm while Heyman and Reeves favoured algorithm appears to have values in the range of 7.2 to 7.4. Thus, the algorithm given in this paper produces results that have not been achieved in the past.

Using the computations for the MFPT matrix *M*, as calculated using the EGTH algorithm, we compute the elements of the stationary distributions as the reciprocal of the diagonal elements. We calculate the following errors for the stationary distribution, both in single and double precision: *Minimum absolute error* $= \min_{1 \leq j \leq m} \left| \pi_j - \sum_{i=1}^{m} \pi_i p_{ij} \right|$, *maximum absolute error* $= \max_{1 \leq j \leq m} \left| \pi_j - \sum_{i=1}^{m} \pi_i p_{ij} \right|$, and the *overall residual error* $= \sum_{j=1}^{m} \left| \pi_j - \sum_{i=1}^{m} \pi_i p_{ij} \right|$, where the $\pi_j$ are given by the calculations. We also compute the overall residual error when the stationary distribution is computed using the standard GTH algorithms. These calculations are given in Table 4 and 5.

Table 4: Errors for the Stationary distributions under double precision

| Test Problem | EGTH Min Abs Error | EGTH Max Abs Error | EGTH Residual Error | GTH Residual Error |
|---|---|---|---|---|
| TP1  | 0 | 1.1102E-16 | 1.4485E-16 | 7.1124E-17 |
| TP2  | 0 | 2.7756E-17 | 7.6328E-17 | 2.0817E-17 |
| TP3  | 0 | 1.3878E-17 | 1.3878E-17 | 1.3878E-17 |
| TP41 | 0 | 2.7756E-17 | 1.1102E-16 | 1.1796E-16 |
| TP42 | 0 | 2.7756E-17 | 8.3267E-17 | 1.0408E-16 |
| TP43 | 0 | 2.7756E-17 | 1.6653E-16 | 1.0408E-16 |
| TP44 | 0 | 2.7756E-17 | 1.1102E-16 | 1.0408E-16 |

Table 5: Errors for the Stationary distributions under single precision

| Test Problem | EGTH Min Abs Error | EGTH Max Abs Error | EGTH Residual Error | GTH Residual Error |
|---|---|---|---|---|
| TP1 | 6.7218E-10 | 2.3568E-08 | 5.4538E-08 | 1.2080E-08 |
| TP2 | 2.1102E-09 | 1.1569E-08 | 5.5893E-08 | 4.9913E-08 |
| TP3 | 8.8180E-15 | 1.4567E-08 | 2.6965E-08 | 2.7865E-08 |
| TP41 | 2.5098E-09 | 2.4648E-08 | 7.3546E-08 | 7.0168E-08 |
| TP42 | 9.4676E-10 | 1.4745E-08 | 6.5571E-08 | 7.0168E-08 |
| TP43 | 1.5393E-09 | 1.16931E-08 | 5.4947E-08 | 6.2717E-08 |
| TP44 | 1.0553E-09 | 1.7522E-08 | 7.9552E-08 | 7.0168E-08 |

As can be expected the errors for computing the stationary distributions using the well established GTH algorithm are very comparable with the EGTH procedure of this paper giving only a marginal reduction but in some isolated cases a slightly improved result.

In order to make comparisons in the case of the stationary distribution calculations that appear in the literature we also compare the errors between performing the calculations for both the EGTH and the original GTH algorithms in double and single precision as follows. Let $\pi_S$ and $\pi_D$ be the stationary distributions as computed under single and double precision. As used in Harrod and Plemmons, ([5]) the *residual error* is, in effect, the residual error computed as above under single precision, i.e, $\left\| \pi_S^T - \pi_S^T P \right\|_1$. The *relative error* is computed as $\sum_{j=1}^{m} \left| \pi_{S,j} - \pi_{D,j} \right|$. We also compute the *minimum absolute error* $\min_{1 \leq j \leq m} \left| \pi_{S,j} - \pi_{D,j} \right|$ and the *maximum absolute error* $\max_{1 \leq j \leq m} \left| \pi_{S,j} - \pi_{D,j} \right|$.

Table 6: Differences between single and double precision computations of the stationary distributions

| Test Problem | ETGH Min Abs Error | EGTH Max Absolute Error | EGTH Relative Error | GTH Relative Error |
|---|---|---|---|---|
| TP1 | 2.3546E-10 | 1.7982E-08 | 4.0117E-08 | 3.8463E-08 |
| TP2 | 7.0444E-10 | 2.8857E-08 | 8.5618E-08 | 5.1491E-08 |
| TP3 | 4.2533E-15 | 1.8365E-08 | 4.8544E-08 | 4.0007E-08 |
| TP41 | 7.0264E-10 | 1.6013E-08 | 6.7861E-08 | 4.5877E-08 |
| TP42 | 7.0264E-10 | 1.1836E-08 | 4.9242E-08 | 4.5877E-08 |
| TP43 | 7.0264E-10 | 1.1836E-08 | 5.5331E-08 | 4.5877E-08 |
| TP44 | 1.0380E-10 | 1.3945E-08 | 6.8623E-08 | 4.5877E-08 |

We make the following observations in respect to each test problem.

**TP1**: The original transition matrix was given as:

$$\begin{bmatrix} .2 & 0 & 0 & .6 & 0 & 0 & 0 & 0 & 0 & .2 \\ 0 & .1 & 0 & 0 & .6 & 0 & .3 & 0 & 0 & 0 \\ 0 & .1 & 0 & 0 & 0 & 0 & 0 & .8 & 0 & .1 \\ 0 & 0 & .6 & 0 & .3 & 0 & 0 & 0 & 0 & .1 \\ 0 & .5 & 0 & 0 & .5 & 0 & 0 & 0 & 0 & 0 \\ 0 & .5 & 0 & 0 & .2 & 0 & 0 & 0 & .3 & 0 \\ 0 & 0 & 0 & 0 & .7 & 0 & .2 & 0 & 0 & .1 \\ .1 & 0 & .9 & 0 & 0 & 0 & 0 & 0 & 0 & 0 \\ 0 & .1 & 0 & 0 & 0 & .8 & 0 & 0 & 0 & .1 \\ 0 & .4 & 0 & 0 & 0 & .4 & 0 & 0 & 0 & .2 \end{bmatrix}$$

Harrod and Plemmons ([5]) gave the exact solution for the solution of the stationary probabilities using some direct methods, however the transition matrix above is not irreducible, and consequently some of the entries of the stationary probability vector should have been zero. Heyman ([6]) commented that the GTH algorithm determines that states 1, 3, 4 and 8 are transient, although this is transparent from an examination of the transition graph. Heyman showed that the GTH algorithm, under single precision, on the above transition matrix produces 6 significant decimal digits (while some alternatives produce only 5) and showed that GTHRE = 4.5 E-08. These were compared with a range of procedures considered by Harrod and Plemmons (1984) that yielded MINRE = 6.9 E-08, MAXRE = 3.7 E-08.

As was done in Heyman and Reeves ([8]), we discard the transient states. With the state space S = {2, 5, 6, 7, 9, 10} we consider the irreducible transition matrix as stated. Using the MatLab single precision our residual error (1.2080E-08) was an improvement over those stated above.

**TP2**: Harrod and Plemmons,[5], stated the exact solution for the stationary distribution to 9 significant figures and showed that the smallest relative error they could achieve was of the order of 9.9 E-07. Heyman, [6], claimed that the GTH algorithm produces 6 significant decimal digits with a residual error of 9.64 E-08. Comparable to the figure of 8.56 E-08 that we have achieved. Under double precision we have been able to achieve 14 significant figures.

**TP3**: Harrod and Plemmons, ([5]), give the exact solution for the stationary distribution to 9 significant figures and using a variety of procedures obtain the smallest residual error of the order of 3.0 E-08. Heyman, using the GTH algorithm, produces 6 significant decimal digits with alternatives giving only 1 or 2. He obtains a residual error of 3.1 E-08 for the GTH algorithm. We have improved this to 14 significant figures, once again with improved accuracy achieving a residual error of 1.4 E-17.

**TP4**: In Harrod and Plemmons [5] the original matrices were not stochastic. Heyman ([6]) corrected this to ensure stochasticity and showed that the stationary distributions of the MCs with these four transition matrices are all the same. He showed that the residual error for the GTH algorithm, under single precision, is 1.38 E-07 for all the four test problems, whereas we achieve accuracy within the range 5.49 E-08 to 7.96E-08.

All in all, the extraction of the stationary distribution as a byproduct of our EGTH algorithm gives comparable accuracy similar to that previously obtained.

In a sequel to this paper we explore some other techniques for computing the MFPTs for these matrices. The results of this paper are required as a benchmark in order to carry out comparisons of accuracy of the alternative procedures.

**Acknowledgement**
The author wishes to acknowledge the assistance of Ms Diane Park who assisted with the MatLab coding as part of her BSc(Hons) dissertation at Auckland University of Technology.

The author although wishes to thank the referees for their comments and recommendations, especially in respect to notation, that has improved the paper.

**References**

[1] Benzi, M. A direct projection method for Markov chains. *Linear Algebra Appl*, **386**, (2004), 27-49.

[2] Bini, D. A., Latouche, G. and Meini B. *Numerical Methods for Structured Markov Chains*, Oxford University Press, New York. (2005).

[3] Dayar, T. and Akar, N. Computing the moments of first passage times to a subset of states in Markov chains. *SIAM J Matrix Anal Appl*, **27**, (2005), 396-412.

[4] Grassman, W.K., Taksar, M.I., and Heyman, D.P. Regenerative analysis and steady state distributions for Markov chains. *Oper Res*, **33**, (1985), 1107-1116.

[5] Harrod, W.J. and Plemmons, R.J. Comparison of some direct methods for computing stationary distributions of Markov chains. *SIAM J Sci Stat Comput*, **5**, (1984), 463-479.

[6] Heyman, D.P. Further comparisons of direct methods for computing stationary distributions of Markov chains. *SIAM J Algebra Discr,* **8**, (1987), 226-232.

[7] Heyman, D.P. and O'Leary, D.P. What is fundamental for Markov chains: First Passage Times, Fundamental matrices, and Group Generalized Inverses, *Computations with Markov Chains*, Chap 10, 151-161, Ed W.J. Stewart, Springer. New York, (1995).

[8] Heyman, D.P. and Reeves, A. Numerical solutions of linear equations arising in Markov chain models. *ORSA J Comp*, **1**, (1989), 52-60.

[9] Hunter, J.J. Generalized inverses and their application to applied probability problems. *Linear Algebra Appl*, **46**, (1982), 157-198.

[10] Hunter, J.J. *Mathematical Techniques of Applied Probability, Volume 2, Discrete Time Models: Techniques and Applications*, Academic, New York. (1983).

[11] Hunter, J.J. Mixing times with applications to Markov chains, *Linear Algebra Appl*, **417**, (2006), 108-123.

[12] Kemeny, J. G. and Snell, J. L. *Finite Markov chains*, Springer- Verlag, New York (1983), (Original version, Princeton University Press, Princeton (1960).)

[13] Kohlas, J. Numerical computation of mean passage times and absorption probabilities in Markov and semi-Markov models. *Zeit Oper Res*, **30,** (1986), 197-207.

[14] Meyer. C. D. Jr. The role of the group generalized inverse in the theory of Markov chains. *SIAM Rev,* **17**, (1975), 443-464.

[15] Sheskin, T.J. A Markov partitioning algorithm for computing steady state probabilities. *Oper Res*, **33**, (1985), 228-235.

[16] Stewart, W. J. *Introduction to the Numerical Solution of Markov chains*. Princeton University Press, Princeton. (1994).

## Appendix 1: MatLab Code for calculations

The code below is an implementation of the EGTH algorithm for the MFPTs and the stationary distribution in the Markov chain setting. Minor modifications can be implemented for the MRP situation.

```
clear all
format long
m=
TM =
e=ones(m,1);
et= ones(1,m);
S=ones(1,m);
E=ones(m,m);
mu=zeros(m,m);
mu(:,m)=1;
PP=TM;
M=zeros(m,m);
P=TM;
   for k=1:m
      for n=m:-1:2
      S(1,n)=sum(PP(n,1:n-1));
      for i=1:n-1
        for j=1:n-1
           PP(i,j)=PP(i,j)+PP(i,n)*PP(n,j)/S(1,n);
        end
         mu(i,n-1)=mu(i,n)+mu(n,n)*PP(i,n)/S(1,n);
      end
      end
M(1,k)=(PP(2,1)*mu(1,2)+PP(1,2)*mu(2,2))/PP(2,1);
      for n=2:m
         mm=0;
         for i=2:n-1
```

```
              mm=mm+PP(n,i)*M(i,k);
          end
           M(n,k)=(mm+mu(n,n))/S(1,n);
        end
        for col=1:m
          for row= 1:m
            P_new(mod(row+m-2,m)+1, mod(col+m-2,m)+1)=P(row,col);
          end
        end
P=P_new;
   PP=P;
   end
   for col=1:m
     for row=1:m
        M_EGTH(mod(row+col-2,m)+1,col)=M(row,col);
     end
   end
M_EGTH
D=diag(diag(M_EGTH));
PI=eye(m)/D;
pit=et*PI
deltaSD=pit-pit*TM;
MINSD=min(abs(deltaSD))
MAXESD=max(abs(deltaSD))
RESD=sum(abs(deltaSD))
MError=M_EGTH-(P*(M_EGTH-D))-E
MinErrorM_EGTH=min(min(abs(MError)))
MaxErrorM_EGTH=max(max(abs(MError)))
REM_EGTH=sum(sum(abs(MError)))
```

# Appendix 2: Mean First Passage Time Matrices and Stationary Probability vectors for the Test Problems

**Test Problem TP1**

*M_GTH (TP1)* = 1.0e+02 *

Columns 1 through 5

```
0.031906040268456   0.016777874480224   1.678333333333333   0.073333333333333   5.641111111111115
0.020000000000000   0.018388937240112   1.698333333333333   0.093333333333333   5.661111111111114
0.023322147651007   0.030226767237211   1.267733333333333   0.096655480984340   3.962777777777780
0.033020134228188   0.017000290107340   1.605000000000000   0.085082774049217   5.567777777777780
0.031073825503356   0.039459433323663   0.263000000000000   0.104407158836689   4.225777777777779
0.024161073825503   0.036002320858718   0.851666666666667   0.097494407158837   4.814444444444447
```

Column 6

```
0.826666666666667
0.846666666666667
0.816403508771930
0.753333333333333
0.745789473684211
0.667228070175439
```

$\pi^T$(TP1) = (0.313420277660917, 0.543805216659655, 0.007888094236432, 0.117532604122844, 0.002366428270930, 0.014987379049222).

**Test Problem TP2**

*M_GTH (TP2)* = 1.0e+04 *

Columns 1 through 5

```
0.001120038405164   0.004000524241398   0.003226904694844   0.239882165430405   0.239916492340981
0.004417409689462   0.001078078341468   0.002304381884289   0.239935793021461   0.239962993464962
0.004262490336574   0.000781374666630   0.002469848581477   0.240060202302811   0.240088526508469
0.301513515251386   0.300259753492525   0.300637382978007   0.000630782736933   0.001031137061417
0.301297906906960   0.300048436220740   0.300423478300436   0.001073299587587   0.000840772722264
1.135651813460504   1.134200327913156   1.134620596021624   1.133546713985028   1.133553779847362
1.135644929057248   1.134193996457640   1.134613896678193   1.133551932417630   1.133559210389771
1.135651812313103   1.134200326857913   1.134620594905067   1.133546714854767   1.133553780752453
```

Columns 6 through 8

```
1.000922840539885   1.000490838593494   1.001389537679506
1.000922873212188   1.000490815512089   1.001389562022972
1.000922957110259   1.000490756645686   1.001389623432740
1.000929150415075   1.000470049184145   1.001438797333724
1.000928904046822   1.000470178394119   1.001438736100158
0.000830664952960   0.000800112187562   0.001636856945469
0.001692132143070   0.000359977354250   0.001819133106368
0.001384723610290   0.000799978835531   0.000982132394897
```

$\pi^T$(TP2) = (0.08928265275450, 0.09275763750513, 0.04048831201636, 0.15853319081983, 0.11893820690418, 0.12038548110605, 0.27779525244927, 0.10181926644467).

**Test Problem TP3**

M_GTH (TP3) = 1.0e+07 *

```
0.000000317241441  0.777777777777778  0.879704977859779  0.341448457627119  0.216494865979381
0.103448448275862  0.511111211111111  0.912915468634687  0.327581118644068  0.123711494845361
0.220689689655172  0.822222222222222  0.000001018450384  0.425269728813559  0.208247432989691
0.220689689655172  0.822222222222222  1.018450383763838  0.000000425269729  0.208247432989691
0.241379379310345  0.666666666666667  0.957195789667897  0.309091000000000  0.000000284536138
```

$\pi^T$ (TP3) = ( 0.31521732963139, 0.00000019565214, 0.09818838658633, 0.23514488152962, 0.35144920660052).

**Test Problem TP4**

M_GTH (TP41) =

Columns 1 through 5

```
 8.045951577090367  12.495886427122791  35.422216573462130  17.601173529910184  11.856207127752848
 3.076923076923078  10.122326177629816  31.602609727164882  13.781566683612947   7.062090618391871
 3.658119658119659   8.886694983354381  26.896466700559230   9.075423657007287   9.567319098457888
 2.743589743589744   9.788992844296484  31.269276393831554  13.448233350279617  10.139541105781630
 2.136752136752137   9.441702062526305  35.260125402474159  17.439082358922217  10.232351462169268
45.673190984578881  58.169077411701679  81.095407558041003  63.274364514489051  57.529398112331720
49.913997627520750  62.409884054643548  85.336214200982880  67.515171157430927  61.770204755273589
50.388493475682083  62.884379902804881  85.810710049144220  67.989667005592253  62.244700603434922
48.935349940688013  61.431236367810811  84.357566514150150  66.536523470598183  60.791557068440852
50.385527876631073  62.881414303753871  85.807744450093210  67.986701406541243  62.241735004383912
```

Columns 6 through 10

```
34.786324786324791  45.246329587753088  49.815688742069113  42.594762594762599  49.286887212516078
37.863247863247864  48.323252664676161  52.892611818992194  45.671685671685680  52.363810289439151
38.444444444444450  48.904449245872748  53.473808400188773  46.252882252882266  52.945006870635730
37.529914529914535  47.989919331342840  52.559278485658872  45.338352338352358  52.030476956105836
36.923076923076927  47.383081724505232  51.952440878821257  44.731514731514736  51.423639349268221
 8.045951577090367  10.460004801428308  15.029363955744323   7.808437808437810  14.500562426191298
 4.240806642941873   7.398677043345709  12.853154003460752   5.064225064225065  10.352225244678761
 4.715302491103203   5.038329848801624  13.870628178910389   5.275135275135275   9.888075629964071
 3.262158956109134   8.105137295210925  13.355356300141578   7.406756406756408  11.309533825027810
 4.712336892052194   6.856509548083042  11.350217607886318   4.810894810894811   9.890976565056040
```

M_GTH (TP42) = 1.0e+03 *

Columns 1 through 5

```
   0.008045951577090   0.012495886427123   0.035422216573462   0.017601173529910   0.011856207127753
   0.003076923076923   0.010122326177630   0.031602609727165   0.013781566683613   0.007062090618392
   0.003658119658120   0.008886694983354   0.026896466700559   0.009075423657007   0.009567319098458
   0.002743589743590   0.009788992844296   0.031269276393832   0.013448233350280   0.010139541105782
   0.002136752136752   0.009441702062526   0.035260125402474   0.017439082358922   0.010232351462169
   4.567319098457888   4.579814984885012   4.602741315031350   4.584920271987799   4.579175305585641
   4.571559905100830   4.584055791527955   4.606982121674293   4.589161078630741   4.583416112228583
   4.572034400948992   4.584530287376116   4.607456617522455   4.589635574478902   4.583890608076744
   4.570581257413997   4.583077143841122   4.606003473987461   4.588182430943907   4.582437464541750
   4.572031435349941   4.584527321777063   4.607453651923404   4.589632608879850   4.583887642477694
```

Columns 6 through 10

```
   3.478632478632480   3.489092483433907   3.493661842588224   3.486440916440916   3.493133041058670
   3.481709401709403   3.492169406510831   3.496738765665147   3.489517839517839   3.496209964135593
   3.482290598290600   3.492750603092027   3.497319962246344   3.490099036099036   3.496791160716790
   3.481376068376070   3.491836073177498   3.496405432331814   3.489184506184507   3.495876630802261
   3.480769230769232   3.491229235570660   3.495798594724976   3.488577668577669   3.495269793195423
   0.008045951577090   0.010460004801428   0.015029363955744   0.007808437808438   0.014500562426191
   0.004240806642942   0.007398677043346   0.012853154003461   0.005064225064225   0.010352225244679
   0.004715302491103   0.005038329848802   0.013870628178910   0.005275135275135   0.009888075629964
   0.003262158956109   0.008105137295211   0.013355356300142   0.007406756406756   0.011309533825028
   0.004712336892052   0.006856509548083   0.011350217607886   0.004810894810895   0.009890976565056
```

M_GTH (TP43)  = 1.0e+05 *

Columns 1 through 5

```
   0.000080459515771   0.000124958864271   0.000354222165735   0.000176011735299   0.000118562071278
   0.000030769230769   0.000101223261776   0.000316026097272   0.000137815666836   0.000070620906184
   0.000036581196581   0.000088866949834   0.000268964667006   0.000090754236570   0.000095673190985
   0.000027435897436   0.000097889928443   0.000312692763938   0.000134482333503   0.000101395411058
   0.000021367521368   0.000094417020625   0.000352601254025   0.000174390823589   0.000102323514622
   4.567319098457888   4.567444057322159   4.567673320623622   4.567495110193187   4.567437660529166
   4.567361506524317   4.567486465388588   4.567715728690051   4.567537518259616   4.567480068595595
   4.567366251482798   4.567491210347070   4.567720473648532   4.567542263218098   4.567484813554076
   4.567351720047449   4.567476678911720   4.567705942213183   4.567527731782747   4.567470282118726
   4.567366221826807   4.567491180691079   4.567720443992542   4.567542233562107   4.567484783898085
```

Columns 6 through 10

```
   3.478632478632479   3.478737078680493   3.478782772272036   3.478710563010563   3.478777484256741
   3.478663247863248   3.478767847911262   3.478813541502805   3.478741332241332   3.478808253487510
   3.478669059829060   3.478773659877075   3.478819353468617   3.478747144207145   3.478814065453322
   3.478659914529916   3.478764514577929   3.478810208169472   3.478737998908000   3.478804920154177
   3.478653846153847   3.478758446201860   3.478804139793404   3.478731930531930   3.478798851778109
   0.000080459515771   0.000104600048014   0.000150293639557   0.000078084378084   0.000145005624262
   0.000042408066429   0.000073986770433   0.000128531540035   0.000050642250642   0.000103522252447
   0.000047153024911   0.000050383298488   0.000138706281789   0.000052751352751   0.000098880756300
   0.000032621589561   0.000810051372952   0.000133553563001   0.000074067564068   0.000113095338250
   0.000047123368921   0.000068565095481   0.000113502176079   0.000048108948109   0.000098909765651
```

M_GTH (TP44) = 1.0e+07 *

Columns 1 through 5

```
0.000000804595158   0.000001249588643   0.000003542221657   0.000001760117353   0.000001185620713
0.000000307692308   0.000001012232618   0.000003160260973   0.000001378156668   0.000000706209062
0.000000365811966   0.000000888669498   0.000002689646670   0.000000907542366   0.000000956731910
0.000000274358974   0.000000978899284   0.000003126927639   0.000001344823335   0.000001013954111
0.000000213675214   0.000000944170206   0.000003526012540   0.000001743908236   0.000001023235146
4.567319098457888   4.567320348046532   4.567322640679545   4.567320858575242   4.567320284078601
4.567319522538552   4.567320772127196   4.567323064760209   4.567321282655906   4.567320708159265
4.567319569988138   4.567320819576780   4.567323112209794   4.567321330105491   4.567320755608850
4.567319424673784   4.567320674262427   4.567322966895441   4.567321184791138   4.567320610294496
4.567319569691578   4.567320819280220   4.567323111913234   4.567321329808930   4.567320755312290
```

Columns 6 through 10

```
3.478632478632480   3.478633524632959   3.478633981568875   3.478633259476260   3.478633928688722
3.478632786324787   3.478633832325267   3.478634289261183   3.478633567168568   3.478634236381030
3.478632844444446   3.478633890444925   3.478634347380841   3.478633625288226   3.478634294500688
3.478632752991454   3.478633798991934   3.478634255927850   3.478633533835235   3.478634203047697
3.478632692307694   3.478633738308173   3.478634195244089   3.478633473151474   3.478634142363936
0.000000804595158   0.000001046000480   0.000001502936396   0.000000780843781   0.000001450056243
0.000000424080664   0.000000739867704   0.000001285315400   0.000000506422506   0.000001035222524
0.000000471530249   0.000000503832985   0.000001387062818   0.000000527513528   0.000000988807563
0.000000326215896   0.000000810513730   0.000001335535630   0.000000740675641   0.000001130953383
0.000000471233689   0.000000685650955   0.000001135021761   0.000000481089481   0.000000989097657
```

$\pi^T$(TP41) = $\pi^T$(TP42) = $\pi^T$(TP43) = $\pi^T$(TP44) =
(0.12428610717064, 0.09879152108435, 0.03717960470916, 0.07435920941833, 0.09772924666409,
0.12428610717064, 0.13515929863426, 0.07209478814524, 0.13501186552967, 0.10110225147362).